\documentclass[10pt]{elsarticle}
\usepackage[cp1251]{inputenc}
\usepackage[english]{babel}
\usepackage{amsmath}
\usepackage{amssymb}
\usepackage{amsfonts}
\usepackage{amsthm}
\usepackage{mathtools}
\usepackage{dsfont}
\usepackage{rotating}
\usepackage{graphicx}
\usepackage{floatflt,epsfig}
\usepackage{lineno,hyperref}
\usepackage{enumerate}
\usepackage{colortbl}
\usepackage{array,tabularx,tabulary,booktabs}
\usepackage{longtable}
\usepackage{multirow}
\usepackage{wrapfig}
\usepackage{subcaption}
\usepackage{pdflscape}
\usepackage[table]{xcolor}
\newcolumntype{^}{>{\currentrowstyle}}

\journal{Arxiv}
\newtheorem{lemma}{Lemma}

\newtheorem{theorem}{Theorem}
\newtheorem{corollary}{Corollary}
\newtheorem{proposition}{Proposition}

\newtheorem{problem}{Problem}

\newcommand{\PG}{\operatorname{PG}}

\begin{document}
\renewcommand{\abstractname}{Abstract}
\renewcommand{\refname}{References}
\renewcommand{\tablename}{Table}
\renewcommand{\arraystretch}{0.9}
\thispagestyle{empty}
\sloppy

\begin{frontmatter}
\title{Divisible design graphs from the symplectic graph} 

\author[01]{Bart De Bruyn}
\ead{Bart.DeBruyn@ugent.be}

\author[02]{Sergey Goryainov}
\ead{sergey.goryainov3@gmail.com}

\author[03]{Willem H.  Haemers}
\ead{haemers@uvt.nl}

\author[04]{Leonid Shalaginov}
\ead{44sh@mail.ru}

\address[01] {Department of Mathematics: Algebra and Geometry, Ghent University, Belgium}
\address[02] {School of Mathematical Sciences, Hebei International Joint Research Center for Mathematics and Interdisciplinary Science, Hebei Key Laboratory of Computational Mathematics and Applications, Hebei Normal University, P.R. China}
\address[03] {Department of Econometrics and Operations Research, Tilburg University, The Netherlands}
\address[04] {Chelyabinsk State University, Russia}

\begin{abstract}
A divisible design graph is a graph whose adjacency matrix is an incidence matrix of a (group) divisible design. Divisible design graphs were introduced in 2011 as a generalization of $(v,k,\lambda)$-graphs. Here we describe four new infinite families that can be obtained from the symplectic strongly regular graph $Sp(2e,q)$ ($q$ odd, $e\geq 2$) by modifying the set of edges. To achieve this we need two kinds of spreads in $PG(2e-1,q)$ with respect to the associated symplectic form: the symplectic spread consisting of totally isotropic subspaces and, when $e=2$, a special spread consisting of lines which are not totally isotropic. Existence of symplectic spreads is known, but the construction of a special spread for every odd prime power $q$ is a major result of this paper. We have included relevant back ground from finite geometry, and when $q=3,5$ and $7$ we worked out all possible special spreads.
\end{abstract}

\begin{keyword}
divisible design graph; symplectic graph; spread; equitable partition; projective space
\vspace{\baselineskip}
\MSC[2020] 	05B25 \sep 05E30 \sep 51A50
\end{keyword}

\end{frontmatter}

\section{Introduction}\label{sec:Intro}
A $k$-regular graph on $v$ vertices is called a \emph{divisible design graph} with parameters $(v,k,\lambda_1,\lambda_2,m,n)$ is its vertex-set can be partitioned into $m$ classes of size $n$ such that any two vertices from the same class have $\lambda_1$ common neighbours and any two vertices from different classes have $\lambda_2$ common neighbours.
Divisible design graphs were introduced in \cite{HKM11} as a bridge between (group) divisible designs and graphs. 
It follows that the adjacency matrix of a divisible design graph is an incidence matrix of a divisible design. 
Thus divisible design graphs provide an interplay for these two areas of combinatorics. They then were subsequently studied in \cite{CH14}, \cite{GHKS19}, \cite{KS21}, \cite{S21}, \cite{PS22}, \cite{CS22}, \cite{P22}, \cite{K22}, and \cite{K23}. In particular, several constructions of divisible design graphs were introduced. We also note that in \cite{PS22} divisible design graphs with at most 39 vertices were enumerated. 

In this paper, for any odd prime power $q$, we construct a new divisible design graph that is based on the symplectic graph $Sp(4,q)$. 
We also show that the complement of the symplectic graph $Sp(2e,q)$ admits equitable partitions that satisfy the requirements of \cite[Construction 4.16]{HKM11}. 
This gives rise to three more infinite families of divisible design graphs (see Theorems \ref{thm:DDG2}, \ref{thm:DDG3} and \ref{thm:DDG4}).
The smallest graphs in these four families have 40 vertices and thus cannot be found in \cite{PS22}. 
Also, these graphs cannot be found in the database of Cayley-Deza graphs with fewer than 60 vertices (see \cite{GS14} and \cite{P}) and thus are not Cayley graphs. 
But we note that some of these graphs can be found in the database of vertex-transitive graphs with fewer than 48 vertices \cite{RH20}. 
Also, the complement of one of these graphs on 40 vertices can be found in \cite[Section 5.6]{BO09} under the name $\Gamma_{8,9}$, and some other graphs of order 40 from these families can be found in~\cite[Section~4.5]{HKM11}. 

Inspired by the ideas from \cite[Section 5.6]{BO09}, we present the following construction, which requires the notion of a special spread in $\operatorname{PG}(3,q)$. A spread $S$ in $\operatorname{PG}(3,q)$, where $q$ is odd, consisting of lines that are not totally isotropic (w.r.t. a given symplectic polarity), having the property that, for any line $\ell \in S$, there exists a uniquely determined line $\ell' \in S$ such that any point of $\ell$ is orthogonal to any point of $\ell'$, is called \emph{special}.
Note that any two paired lines $\ell,\ell'$ from a special spread $S$ induce a complete bipartite graph $K_{q+1,q+1}$ in the symplectic graph $Sp(4,q)$; moreover, there exists a partition (associated with $S$) of the vertex set of $Sp(4,q)$ into such $K_{q+1,q+1}$'s.

\begin{theorem}\label{thm:DDG1}
Consider $Sp(4,q)$ with $q$ odd and a special spread $S$.
Let $\Gamma_{q,S}$ be the complement of the graph obtained from $Sp(4,q)$ by removing all the edges of each $K_{q+1,q+1}$ in the partition of the vertex set $V(Sp(4,q))$ into $K_{q+1,q+1}$'s associated with the spread $S$. Then the graph $\Gamma_{q,S}$ is a divisible design graph with parameters 
$$
((q^2+1)(q+1),q^3+q+1,q^3-q^2+q+1,q^3-q^2+2q,q^2+1,q+1).
$$
\end{theorem}

With a special spread we can make other divisible design graphs.

\begin{theorem}\label{thm:DDG2}
Consider $Sp(4,q)$ with $q$ odd and a special spread $S$.
Partition the vertices of $Sp(4,q)$ into two parts $V_1$ and $V_2$ of equal size, such that, for every subgraph $K_{q+1,q+1}$ associated with $S$, one part is in $V_1$ and the other part is in $V_2$.
Let $\Gamma'_{q,S}$ be the graph obtained from $Sp(4,q)$ by replacing the subgraphs induced by $V_1$ and $V_2$ by their complements.
Then $\Gamma'_{q,S}$ is a divisible design graph with parameters
$$
((q^2+1)(q+1),(q^3+q^2+3q+1)/2,(q^3-q^2+3q+1)/2,q^2+q,2,(q^2+1)(q+1)/2).
$$
\end{theorem}

Note that there are exponentially many choices for such a partition into $V_1$ and $V_2$. Therefore the above construction gives many nonisomorphic divisible design graphs.

The following theorem is the core of the previous constructions since it shows the existence (in a constructive manner) of a special spread for each odd prime power $q$.

\begin{theorem}\label{thm:spspr}
Given an odd prime power $q$, there exists at least one special spread in $\operatorname{PG}(3,q)$.
\end{theorem}

With the aid of the computer algebra systems GAP \cite{GAP} and SageMath \cite{SAGE}, we have verified in a computational way the following result, see Section \ref{sec:SpecialSpreads1.5}.

\begin{theorem}\label{thm:computer}
For $q = 3,5$ and $7$, the projective space $\operatorname{PG}(3,q)$ has exactly $1,2$ and $14$ pairwise non-equivalent special spreads.
\end{theorem}

Here, two special spreads are called {\em equivalent} if there exists an automorphism of $Sp(4,q)$ mapping one of them to the other. We have also verified that the non-equivalent special spreads give (in Theorem \ref{thm:DDG1}, but not necessarily in Theorem \ref{thm:DDG2}) non-isomorphic divisible design graphs. We generalise this phenomenon in the following theorem.

\begin{theorem}\label{thm:2spspr}
Let $q$ be an odd prime power, and $S_1$ and $S_2$ be two non-equivalent special spreads in $\operatorname{PG}(3,q)$. Then the graphs $\Gamma_{q,S_1}$ and $\Gamma_{q,S_2}$ are not isomorphic.
\end{theorem}

It is known that $PG(3,q)$ also has a so called \emph{symplectic spread} consisting of totally isotropic lines partitioning the point set.
This spread corresponds to a partition of the vertices of $Sp(4,q)$ into $q^2+1$ cliques of order $q+1$.
More generally, $Sp(2e,q)$ with $e \geq 2$ has a symplectic spread consisting of totally isotropic subspaces, which
corresponds to a partition of the vertices of $Sp(2e,q)$ into $q^e+1$ cliques of order $(q^e-1)/(q-1)$;
see e.g. \cite{D77}.
With these spreads we can make two more families of divisible design graphs.

\begin{theorem}\label{thm:DDG3}
Consider $Sp(2e,q)$ with $e \geq 2$ and a symplectic spread $R$.
Let $\Gamma_{q,e,R}$ be the graph obtained from $Sp(2e,q)$ by removing the edges of the cliques in the spread. Then $\Gamma_{q,e,R}$ is a divisible design graph with parameters
$$
\left(
\frac{q^{2e}-1}{q-1},\ q^e\frac{q^{e-1}}{q-1},
\ q^e\frac{q^{e-2}-1}{q-1},\ \frac{(q^{e-1}-1)^2}{q-1},
\ q^e+1,\ \frac{q^e-1}{q-1}
\right).
$$
\end{theorem}

\begin{theorem}\label{thm:DDG4}
Consider $Sp(2e,q)$ with $e \geq 2$ and $q$ odd, and a symplectic spread $R$.
Partition the vertices of $Sp(2e,q)$ into two parts $V_1$ and $V_2$ of equal size, such that each part contains $(q^e+1)/2$ cliques of the spread.
Let $\Gamma'_{q,e,R}$ be the graph obtained from $Sp(2e,q)$ by replacing the subgraphs induced by $V_1$ and $V_2$ by their complements.
Then $\Gamma'_{q,e,R}$ is a divisible design graph with parameters
$$
\left(
v=\frac{q^{2e}-1}{q-1},\ \frac{v}{2}-q^{e-1},\ \frac{v}{2} -q^{2e-2}-q^{e-1},\ q^{2e-2}-q^{e-1},\ 2,\ \frac{v}{2}
\right).
$$
\end{theorem}

Also in Theorem~\ref{thm:DDG4} we obtain many non-isomorphic divisible design graphs with these parameters, because there are exponentially many choices for the partition into $V_1$ and $V_2$.
However, in Theorem~\ref{thm:DDG3} we find just one
divisible design graph for a given spread.
\\

The weight-distribution bound is a lower bound on the cardinality of support of an eigenfunction of a distance-regular graph corresponding to a non-principal eigenvalue (see \cite{KMP16} and \cite{SV21}). Recently, the tightness of the weight-distribution bound was shown for the smallest eigenvalue of a class of generalised Paley graphs of square order \cite{GSY23}, affine polar graphs \cite{GY23} and the block graphs of geometric Steiner systems \cite{GP23}. 
In connection with the complete bipartite graphs $K_{q+1,q+1}$ occurring in the notion of a special spread, we also show the following.

\begin{theorem}\label{thm:WD}
Let $q$ be a prime power ($q$ may be even or odd). 
The weight-distribution is tight for the negative non-principal eigenvalue $-(q+1)$ of $Sp(4,q)$.
\end{theorem}

As a byproduct of the proof of Theorem \ref{thm:WD}, we also find all optimal eigenfunctions for the eigenvalue $-(q+1)$.

This paper is organised as follows. In Section \ref{sec:prelim}, we give necessary definitions and preliminary results. In Section \ref{sec:FromSpecialSpread}, we prove Theorems \ref{thm:DDG1}, \ref{thm:DDG2}, \ref{thm:DDG3}, \ref{thm:DDG4} and \ref{thm:WD}.
In Section \ref{sec:SpecialSpreads1}, we prove Theorem \ref{thm:spspr}. In section \ref{sec:SpecialSpreads2}, we prove Theorem \ref{thm:2spspr}. In Section \ref{sec:SpecialSpreads1.5}, we prove Theorem \ref{thm:computer}. In Section \ref{sec:Discussion}, we discuss the obtained results and formulate related open problems. 

\section{Preliminaries}\label{sec:prelim}
In this section, we give necessary definitions and preliminary results.

\subsection{Equitable partitions}
Let $G$ be a $k$-regular graph with the vertex set $V(G)$. 
Let $\Pi := (V_1,\ldots,V_t)$ be a partition of $V(G)$ into $t$ parts ($t$-partition). 
The partition $\Pi$ is said to be an \emph{equitable} $t$-partition if for any $i,j \in \{1,\ldots,t\}$ there is a constant $p_{ij}$ such that any vertex from the part $V_i$ is adjacent to precisely $p_{ij}$ vertices from the part $V_j$.
The square matrix $P_\Pi:=(p_{ij})_{i,j = 1}^t$ is called the {\em quotient matrix} of the equitable $t$-partition $\Pi$.
Since all row sums of the adjacency matrix $A$ of $G$ and the quotient matrix $P_\Pi$ are equal to $k$, both matrices have eigenvalue $k$. 
Moreover, it is well-known (see for example \cite[Theorem 9.3.3]{GR01}),
that every eigenvalue of $P$ is an eigenvalue of $A$.
The eigenvalue $k$ of $A$ is called \emph{principal}.
An eigenvalue $\theta$ of $A$ is called \emph{non-principal} if $\theta \ne k$. If $\Pi$ is an equitable $2$-partition,
then precisely one non-principal eigenvalue $\theta$ is an eigenvalue of the quotient matrix $P_\Pi$.
In this case we say that the equitable $2$-partition $\Pi$ is $\theta$-\emph{equitable}.

A divisible design graph with parameters $(v,k,\lambda_1,\lambda_2,m,n)$
is called \emph{proper} if $m>1$, $n>1$ and $\lambda_1\neq \lambda_2$.
Otherwise it is called \emph{improper}.
The $m$-partition of its vertex set is then called \emph{canonical}.
It was shown in \cite[Theorem 3.1]{HKM11} that the canonical partition of a proper divisible design graph is always equitable.

\subsection{Strongly regular graphs}

A $k$-regular graph on $v$ vertices is called \emph{strongly regular} with parameters $(v,k,\lambda,\mu)$ if any two adjacent vertices have $\lambda$ common neighbours and any two distinct non-adjacent vertices have $\mu$ common neighbours. 
If $G$ is a strongly regular graph, then its complement $\overline{G}$ is also strongly regular.
The parameters of $\overline{G}$ are $(v,\overline{k},\overline{\lambda},\overline{\mu})=(v,v-k-1,v-2k+\mu-2,v+\lambda-2k)$. 
A strongly regular graph $G$ is \emph{primitive} if both $G$ and $\overline{G}$
are connected. 
If $G$ is not primitive, we call it \emph{imprimitive}. The imprimitive strongly regular graphs are exactly the disjoint unions of complete graphs of the same size and
their complements, namely, the complete multipartite graphs with multipartite parts of equal size. So imprimitive strongly regular graphs are (rather trivial) divisible design graphs.

\begin{lemma}[{\cite[Theorem 5.2.1]{GM15}}]\label{EigenvaluesSRG}
If $G$ is a primitive strongly regular graph with parameters $(v,k,\lambda,\mu)$ and  
$$
\Delta:=\sqrt{(\lambda-\mu)^2+4(k-\mu)},
$$
then $G$ has exactly three eigenvalues 
$$
k,~~~r = \frac{\lambda-\mu+\Delta}{2},~~~s = \frac{\lambda-\mu-\Delta}{2},
$$
with respective multiplicities
$$
m_k = 1,~~~m_{r} = -\frac{(v-1)s+k}{r-s},~~~
m_{s} = \frac{(v-1)r+k}{r-s}.
$$
\end{lemma}
Note that $s < 0 < r$ holds. 
For the primitive strongly regular graph $G$ from the above lemma, the matrix
$$
\left(
  \begin{array}{c|cc}
    1 & k & v-1-k\\
    m_{r} & r & -1-r \\
    m_{s} & s & -1-s\\
  \end{array}
\right)
$$
is called the \emph{modified matrix of eigenvalues}. 
The first column gives the dimensions of the eigenspaces (i.e., the multiplicities of the eigenvalues); the second column contains the eigenvalues of $G$, and the third gives the eigenvalues of its complement $\overline{G}$.

If the parameters of a strongly regular graph $G$ satisfy $\lambda=\mu$, then $G$ is called a \emph{$(v,k,\lambda)$-graph}
(because the adjacency matrix is also an incidence matrix of a symmetric $(v,k,\lambda)$-design).
If $\lambda=\mu-2$, then $\overline{G}$, the complement of $G$,
is a $(v,\overline{k},\overline{\lambda})$-graph.
Note that an improper divisible design graph is a $(v,k,\lambda)$-graph with $\lambda=\lambda_1=\lambda_2$.

Delsarte proved \cite{D73} that the clique number of a strongly regular graph $G$ is at most $1-\frac{k}{s}$. 
A clique in a strongly regular graph whose size attains this
bound is called a \emph{Delsarte clique}.
A partition of the vertex set of a strongly regular graph $G$ into Delsarte cliques
is called a \emph{spread}, see~\cite{HT96}.
Such a pertition is always equitable.
For the complement of $G$, a spread corresponds to an equitable partition into cocliques.
This partition gives a coloring of $\overline{G}$ for which the number of colors 
meets Hoffman's lower bound for the chromatic number.
Therefore it is called a \emph{Hoffman coloring}.

\subsection{The weight-distribution bound for strongly regular graphs} \label{wdb}

Let $\theta$ be an eigenvalue of a graph $G$. A real-valued function $f$ on the vertex set of $\Gamma$ is called an \emph{eigenfunction}
of the graph $G$ corresponding to the eigenvalue $\theta$ (or a \emph{$\theta$-eigenfunction} of $G$), if $f \not \equiv 0$ and
for any vertex $u$ in $G$ the condition
\begin{equation}\label{LocalCondition}
\theta\cdot f(u)=\sum_{\substack{w\in{G(u)}}}f(w)
\end{equation}
holds, where $G(u)$ is the set of neighbours of the vertex $u$.
Although eigenfunctions of graphs receive less attention of researchers in contrast to their eigenvalues, there are still tons of related literature. We refer to the recent survey \cite{SV21} for a summary of results on the problem of finding the minimum cardinality of support of eigenfunctions of graphs and characterising the optimal eigenfunctions.

The following lemma gives lower bounds (the so-called {\em weight-distribution bounds}) for the number of non-zeroes (i.e., the cardinality of the support) for an eigenfunction of a strongly regular graph. In fact, this is a special case of a more general result for distance-regular graphs \cite[Section 2.4]{KMP16}.

\begin{lemma}[{\cite[Corollary 3]{GP23}}]\label{WDBsrg}
Let $G$ be a primitive strongly regular graph with non-principal eigenvalues $s<0<r$. Then an eigenfunction of $G$ corresponding to the eigenvalue $r$ has at least $2(r+1)$ non-zeroes, and an eigenfunction corresponding to the eigenvalue $s$ has at least $-2s$ non-zeroes.
\end{lemma}

The following lemma gives a combinatorial interpretation of the tightness of the weight-distribution bound in terms of special induced subgraphs. With isolated cliques $T_0$ and $T_1$, we mean in this lemma that there are no edges between vertices of $T_0$ and vertices of $T_1$.

\begin{lemma}[{\cite[Lemma 3]{GP23}}]\label{OptimalEigenfunctionsSRG}
Let $G$ be a primitive strongly regular graph with non-principal eigenvalues $s,r$, where $s < 0 < r$. Then the following statements hold.\\
{\rm (1)} For an $s$-eigenfunction $f$, if the cardinality of support of $f$ meets the weight-distribution bound, then there exists an induced complete bipartite subgraph in $G$ with parts $T_0$ and $T_1$ of size $-s$. Moreover, up to multiplication by a constant, $f$ has value $1$ on the vertices of $T_0$ and value $-1$ on the vertices of $T_1$. \\
{\rm (2)} For an $r$-eigenfunction $f$, if the cardinality of support of $f$ meets the weight-distribution bound, then there exists an induced pair of isolated cliques $T_0$ and $T_1$ in $G$ of size $-\overline{s} = -(-1-r) = 1+r$. Moreover, up to multiplication by a constant, $f$ has value 1 on the vertices of $T_0$ and value $-1$ on the vertices of $T_1$.\\
{\rm (3)} If $G$ has Delsarte cliques and each edge of $G$ lies in a constant number of Delsarte cliques (for example, $G$ is an edge-transitive strongly regular graph with Delsarte cliques), then any copy (as an induced subgraph) of the complete bipartite graph with parts of size $-s$ in $G$ gives rise to an eigenfunction of $G$ whose cardinality of support meets the weight-distribution bound and which is of the form given in item {\rm(1)}.\\
{\rm (4)} If the complement of $G$ has Delsarte cliques and each edge of the complement of $G$ lies in a constant number of Delsarte cliques (for example, $G$ is a coedge-transitive strongly regular graph whose complement has Delsarte cliques), then any copy (as an induced subgraph) of a pair of isolated cliques of size $1+r$ in $G$ gives rise to an eigenfunction of $G$ whose cardinality of support meets the weight-distribution bound and which is of the form given in item {\rm(2)}.
\end{lemma}

Thus, in view of Lemma \ref{OptimalEigenfunctionsSRG}, to show the tightness of the weight-distribution bound for non-principal eigenvalues, it suffices to find a special induced subgraph (a pair of isolated cliques $T_0$ and $T_1$ or a complete bipartite graph with parts $T_0$ and $T_1$) and show that each vertex outside of $T_0 \cup T_1$ has the same number of neighbours in $T_0$ and $T_1$. Moreover, if we can verify the condition from Lemma \ref{OptimalEigenfunctionsSRG}(3) or \ref{OptimalEigenfunctionsSRG}(4), then it suffices to find a required special induced subgraph.

\subsection{Notions from projective geometry}\label{geometry}

The aim of this subsection is to recall some basic facts about projective geometry, quadrics and polarities. More background information on these topics can be found in the monographs \cite{Hi} and \cite{Hi-Th}.

Let $\mathrm{PG}(n,q)$ be the projective space of dimension $n$ over the finite field $\mathbb{F}_q$. The points of $\mathrm{PG}(n,q)$ are the $1$-dimensional subspaces of an $(n+1)$-dimensional vector space $V$ over $\mathbb{F}_q$. After fixing a basis $(\bar e_1,\bar e_2,\ldots,\bar e_{n+1})$ of $V$, we can denote the point $p=\langle X_1 \bar e_1 + X_2 \bar e_2 + \cdots + X_{n+1} \bar e_{n+1} \rangle$ of $\mathrm{PG}(n,q)$ by its so-called {\em homogeneous coordinates} $(X_1,X_2,\ldots,X_{n+1})$. With every $(i+1)$-dimensional subspace $W$ of $V$, there is associated an $i$-dimensional subspace of $\mathrm{PG}(n,q)$ which is the set of all points that are contained in $W$ (as 1-dimensional subspace). Such a subspace of $\mathrm{PG}(n,q)$ is called a {\em line}, {\em plane}, {\em solid} or {\em hyperplane} depending on whether $i$ is equal to $1$, $2$, $3$ or $n-1$. 

A {\em quadric} $Q$ of $\mathrm{PG}(n,q)$ is a set of points of $\mathrm{PG}(n,q)$ whose homogeneous coordinates $(X_1,X_2,\ldots,X_{n+1})$ (with respect to a basis of $V$) satisfy an equation of the form
\[  \sum_{i,j=1}^{n+1} a_{ij} X_i X_j =0,  \]
where the $a_{ij}$'s are given elements of $\mathbb{F}_q$. A point $x$ of $Q$ is called {\em singular} if $xy \subseteq Q$ for any $y \in Q \setminus \{ x \}$. The quadric $Q$ is called {\em singular} if it has singular points; otherwise it is called {\em nonsingular}. 

If $Q$ is a quadric of $\mathrm{PG}(n,q)$, then a line of $\mathrm{PG}(n,q)$ intersects $Q$ in either 0, 1, 2 or $q+1$ points. Lines intersecting $Q$ in either 1 or $q+1$ points are called {\em tangent lines}, lines intersecting $Q$ in two points are called {\em secant lines} and lines disjoint from $Q$ are called {\em external lines}. A point $x$ of $Q$ is thus singular if and only if all lines through $x$ are tangent lines. If $Q$ is nonsingular, then for every point $x \in Q$, there is a hyperplane $T_x$ of $\mathrm{PG}(n,q)$ through $x$ such that the tangent lines through $x$ are precisely the lines of $T_x$ through $x$. The hyperplane $T_x$ is called the {\em tangent hyperplane at the point} $x$.

In $\mathrm{PG}(n,q)$ with $n \geq 2$ even, there is up to projective equivalence only one nonsingular quadric. The quadric is called a {\em parabolic} quadric, denoted by $Q(n,q)$, and has equation $X_1X_2+X_3X_4+\cdots+X_{n-1}X_n+X_{n+1}^2=0$ with respect to a certain basis of $V$. The quadric $Q(2,q)$ of $\mathrm{PG}(2,q)$ is usually called an {\em irreducible conic}. 

In $\mathrm{PG}(n,q)$ with $n \geq 1$ odd, there are up to projective equivalence two nonsingular quadrics. One of them is called a {\em hyperbolic quadric}, is denoted by $Q^+(n,q)$, and has equation $X_1X_2+\cdots+X_n X_{n+1}=0$ with respect to a certain basis of $V$. The other is called an {\em elliptic quadric}, is denoted by $Q^-(n,q)$, and has equation $X_1X_2+X_3X_4+\cdots+X_{n-2}X_{n-1} + X_n^2 + a X_n X_{n+1} + b X_{n+1}^2 =0$ with respect to a certain basis of $V$. Here, $a$ and $b$ are given elements of $\mathbb{F}_q$ such that the polynomial $X^2+aX+b$ of $\mathbb{F}_q[X]$ is irreducible. The hyperbolic quadric $Q^+(5,q)$ is also called the {\em Klein quadric}.

An {\em anti-automorphism} of $\mathrm{PG}(n,q)$ is a permutation of the set of subspaces of $\mathrm{PG}(n,q)$ reversing the inclusion of subspaces. Such an anti-automorphism is called a {\em polarity} if it has order 2. In $\mathrm{PG}(n,q)$ with $n$ odd there exist polarities $\zeta$ having the property that $x \in x^\zeta$ for every point $x$. Such polarities are called {\em symplectic}. 

With every nonsingular quadric $Q$ of $\mathrm{PG}(n,q)$, $q$ odd, there is associated another polarity $\zeta$. A point $p$ of $\mathrm{PG}(n,q)$ belongs to $Q$ if and only if $p \in p^\zeta$, in which case $p^\zeta$ coincides with the tangent hyperplane at the point $p$. The nontangent hyperplanes to the quadric $Q$ are precisely the hyperplanes $x^\zeta$ for points $x$ not belonging to $Q$. Any such nontangent hyperplane $\Pi=x^\zeta$ intersects $\Pi$ in a nonsingular quadric $Q'$ of $\Pi$, and the tangent lines through $x$ are precisely the lines connecting $x$ with a point of $Q'$. Polarities which are associated with a quadric $Q$ of $\mathrm{PG}(n,q)$, $q$ odd, in the above way are called {\em orthogonal}.

Let $Q(2,q)$ be an irreducible conic in $\mathrm{PG}(2,q)$, $q$ odd, with associated polarity $\zeta$. Let $x$ be point of $\mathrm{PG}(2,q)$ not contained in $Q(2,q)$. The line $x^\zeta$ of $\mathrm{PG}(2,q)$ then intersects $Q(2,q)$ in either 0 or 2 points. If the former case, the point $x$ is called {\em interior} with respect to $Q(2,q)$. In the latter case, $x$ is called {\em exterior} with respect to $Q(2,q)$. Through an exterior point, there are 2 tangent lines, $\frac{q-1}{2}$ lines intersecting $Q(2,q)$ in two points and $\frac{q-1}{2}$ lines disjoint from $Q(2,q)$. Through an interior point, there are no tangent lines, $\frac{q+1}{2}$ lines intersecting $Q(2,q)$ in two points and $\frac{q+1}{2}$ lines disjoint from $Q(2,q)$. There are $\frac{q(q+1)}{2}$ exterior points and $\frac{q(q-1)}{2}$ interior points.

The {\em Klein correspondence} $\kappa$ is a certain nice bijective map between the set of lines of $\mathrm{PG}(3,q)$ and the set of points of the Klein quadric. Line pencils, being sets of lines of $\pi$ through $x$ for incident point-plane pairs $(x,\pi)$, are mapped by $\kappa$ to lines of $Q^+(5,q)$. For any point $x$ of $\mathrm{PG}(3,q)$, $\kappa$ maps the set of lines of $\mathrm{PG}(3,q)$ through $x$ to a plane of $Q^+(5,q)$, called a {\em Latin plane}, and for any plane $\pi$ of $\mathrm{PG}(3,q)$, $\kappa$ maps the set of lines contained in $\pi$ to a plane of $Q^+(5,q)$, called a {\em Greek plane}. Every plane of $Q^+(5,q)$ is either a Latin or a Greek plane. An {\em ovoid} of $Q^+(5,q)$ is a set of points intersecting each plane of $Q^+(5,q)$ in a singleton. The standard examples of ovoids of $Q^+(5,q)$ are given by the elliptic quadrics $Q^-(3,q) \subseteq Q^+(5,q)$ obtained by intersecting $Q^+(5,q)$ with a suitable solid of $\mathrm{PG}(5,q)$. If $O$ is an ovoid of $Q^+(5,q)$, then $\kappa^{-1}(O)$ is a {\em line-spread} of $\mathrm{PG}(3,q)$, being a set of lines of $\mathrm{PG}(3,q)$ partitioning its point set. 

For a given symplectic polarity $\zeta$ of $\mathrm{PG}(3,q)$, the lines $L$ of $\mathrm{PG}(3,q)$ satisfying $L^\zeta=L$ are called {\em totally isotropic} or {\em symplectic}. Any other line of $\mathrm{PG}(3,q)$ is then called {\em hyperbolic}. The Klein correspondence $\kappa$ maps the set of totally isotropic lines to a parabolic quadric $Q(4,q) \subseteq Q^+(5,q)$ obtained by intersecting $Q^+(5,q)$ with a suitable nontangent hyperplane of $\mathrm{PG}(5,q)$ with respect to $Q^+(5,q)$.

The point-line geometry formed by the points and totally isotropic lines of $\mathrm{PG}(3,q)$ is a {\em generalized quadrangle}, meaning that for every non-incident point-line pair $(x,L)$, there exists a unique point on $L$ collinear with $x$. This generalized quadrangle is denoted by $W(q)$ and called a {\em symplectic generalized quadrangle}. The graph whose vertices are the points of $W(q)$, with two distinct points being adjacent whenever they are collinear in $W(q)$, is denoted by $Sp(4,q)$. 

If $K$ is a hyperbolic line, then also $K^\zeta$ is a hyperbolic line. The line $K^\zeta$ is disjoint from $K$ and any line connecting a point of $K$ with a point of $K^\zeta$ is totally isotropic. The hyperbolic lines $K$ and $K^\zeta$ are called {\em orthogonal}. We denote by $\mathcal{U}_q$ the set of all point sets that arise as the union of two orthogonal hyperbolic lines. The point sets of $Sp(4,q)$ on which the induced subgraphs are isomorphic to the complete bipartite graph $K_{q+1,q+1}$ are precisely the elements of $\mathcal{U}_q$.

More generally, for every $e \in \mathbb{N} \setminus \{ 0,1 \}$ and every prime power $q$, we can define a graph $Sp(2e,q)$ whose vertices are the points of $\mathrm{PG}(2e-1,q)$, where two distinct  points $x$ and $y$ are adjacent whenever $y \in x^\zeta$, where $\zeta$ is a given symplectic polarity of $\mathrm{PG}(2e-1,q)$. The graph $Sp(2e,q)$ is called a {\em symplectic graph}. A subspace $\pi$ of $\mathrm{PG}(2e-1,q)$ is called {\em totally isotropic with respect to} $\zeta$ if $\pi \subseteq \pi^\zeta$. The following is known about the symplectic graph $Sp(2e,q)$.

\begin{lemma}[{\cite[Section 2.5]{BV22}}]\label{SpGraphs}
The graph $Sp(2e,q)$ is a rank 3 (in particular, arc-transitive) strongly regular graph with parameters
\[ v = \frac{q^{2e}-1}{q-1},\ k=\frac{q(q^{2e-2}-1)}{q-1},\ \lambda=\frac{q^2(q^{2e-4}-1)}{q-1}+q-1,\
\mu=\frac{k}{q} = \lambda+2, \]
and eigenvalues $r = q^{e-1}-1$, $s = -q^{e-1}-1.$
\end{lemma}

Since $\lambda=\mu-2$ for $Sp(2e,q)$, the complement  $\overline{Sp(2e,q)}$ is
a $(v,\overline{k},\overline{\lambda})$-graph with $\overline{k}=q^{2e-1}$ and $\overline{\lambda}=(q-1)q^{2e-2}$.

\section{Proofs of Theorems \ref{thm:DDG1},~\ref{thm:DDG2},~\ref{thm:DDG3},~\ref{thm:DDG4} and \ref{thm:WD}}\label{sec:FromSpecialSpread}

\begin{lemma} \label{lem:PairedLines}
Let $X$ be a set of vertices of $Sp(4,q)$ on which the induced subgraph is a complete bipartite graph isomorphic to $K_{q+1,q+1}$, i.e. $X$ is of the form $l \cup l'$ where $l$ and $l'$ are two orthogonal hyperbolic lines with respect to the symplectic polarity $\zeta$. Then every vertex $x$ of $Sp(4,q)$ not contained in $X$ is adjacent to precisely two vertices of $X$, one in $l$ and another in $l'$.
\end{lemma}
\begin{proof}
The vertices adjacent to $x$ are exactly the vertices in $x^\zeta \setminus \{ x \}$. As $x \not\in l \cup l'$, the plane $x^\zeta$ contains neither $l$ nor $l'$ and so intersects these lines in unique points.
\end{proof}

\medskip Now let us prove Theorem \ref{thm:DDG1}.

Let $Y_{q,S}$ the graph obtained from $Sp(4,q)$ by removing all the edges of each $K_{q+1,q+1}$ in the partition of the vertex set $V(Sp(4,q))$ into $K_{q+1,q+1}$'s naturally associated with the spread $S$. In what follows, we will implicitly make use of Lemma \ref{lem:PairedLines}. It is then easy to see that $Y_{q,S}$ is regular with degree $k' = q^2-1$.

Let us calculate its intersection numbers in $Y_{q,S}$. Let $x$ and $y$ be two arbitrary vertices. Consider the following four cases that occur.

\textbf{Case 1.} $x$ and $y$ belong to the same part of a former $K_{q+1,q+1}$.
Then $x$ and $y$ are not adjacent in $Y_{q,S}$, have 0 common neighbours in $Y_{q,S}$ and $v-2k'-2 = q^3-q^2+q+1$
common neighbours in $\Gamma_{q,S}$ (the complement of $Y_{q,S}$). 

\textbf{Case 2.} $x$ and $y$ belong to different parts of a former $K_{q+1,q+1}$.
Then $x$ and $y$ are not adjacent in $Y_{q,S}$, have $q-1$ common neighbours in $Y_{q,S}$ and $v-2k'-2+q-1 = q^3-q^2+2q$
common neighbours in $\Gamma_{q,S}$. 

\textbf{Case 3.} $x$ and $y$ belong to different copies of $K_{q+1,q+1}$, $x$ and $y$ are adjacent in $Y_{q,S}$ (as well as in the original graph $Sp(4,q)$).
Then $x$ and $y$ lose two common neighbours after the removing of edges, 
have $q-3$ common neighbours in $Y_{q,S}$ and 
$v-2-2(k'-1)+q-3 = q^3-q^2+2q$ common neighbours in $\Gamma_{q,S}$.

\textbf{Case 4.} $x$ and $y$ belong to different copies of $K_{q+1,q+1}$, $x$ and $y$ are not adjacent in $Y_{q,S}$ (as well as in the original graph $Sp(4,q)$).
Then $x$ and $y$ lose two common neighbours after the removing of edges, 
have $q-1$ common neighbours in $Y_{q,S}$ and $v-2-2k'+q-1 = q^3-q^2+2q$ common neighbours in $\Gamma_{q,S}$.

Thus, $\Gamma_{q,S}$ is a divisible design graph with classes given by parts of the subgraphs $K_{q+1,q+1}$'s involved in the partition associated with the spread $S$. $\square$

\medskip The proofs of Theorems \ref{thm:DDG2}, \ref{thm:DDG3} and \ref{thm:DDG4} are consequences of the following result (Construction 4.16) from \cite{HKM11}.

The \emph{partial complement} of a partitioned graph is obtained by taking the complement only with respect to the edges and non-edges between different classes (and thus leave the subgraphs induced by each class unchanged).

\begin{theorem}\label{partial complement}
Suppose $G$ is a $(v,k,\lambda)$-graph.
If $G$ has an equitable 2-partition, or a partition corresponding to a Hoffman coloring,
then the partial complement is a divisible design graph.
\end{theorem}

We know from Sections~\ref{sec:Intro} and \ref{sec:prelim} that $Sp(2e,q)$ has a spread and that its complement $\overline{Sp(2e,q)}$ has a Hoffman coloring. We observed already that $\overline{Sp(2e,q)}$ is a $(v,\overline{k},\overline{\lambda})$-graph. So Theorem~\ref{partial complement} applies and we find the divisible design graph of Theorem~\ref{thm:DDG3}.

Now let us prove Theorem \ref{thm:DDG4}. Because a spread of $Sp(2e,q)$ is an equitable partition, the partition into $V_1$ and $V_2$ is an equitable 2-partition of $Sp(2e,q)$ as well as its complement $\overline{Sp(2e,q)}$. Now Theorem~\ref{partial complement} implies that the partial complement is a divisible design graph and the parameters given in Theorem~\ref{thm:DDG4} readily follow.

From Lemma~\ref{lem:PairedLines} it follows that a special spread in $Sp(4,q)$ gives an equitable $(q^2+1)$-partition. Therefore the partition into $V_1$ and $V_2$ is an equitable 2-partition of $Sp(4,q)$ and its complement $\overline{Sp(4,q)}$. Now Theorem~\ref{partial complement} implies Theorem~\ref{thm:DDG2}.

\medskip Finally, let us prove Theorem \ref{thm:WD}.

Note that $Sp(4,q)$ is defined on the set of points of $PG(3,q)$, where each line has exactly $q+1$ points. Recall that $Sp(4,q)$ has Delsarte cliques. By Lemma \ref{OptimalEigenfunctionsSRG}(1), if the cardinality of support of an $s$-eigenfunction equals the weight-distribution bound, then the support of this eigenfunction induces a complete bipartite graph with parts of size $-s$. In view of Lemma \ref{OptimalEigenfunctionsSRG}(3) and Lemma \ref{SpGraphs}, any copy of $K_{q+1,q+1}$ as an induced subgraph gives an eigenfunction whose cardinality of support meets the weight-distribution bound. Such subgraphs exist: they are precisely the sets $\mathcal{U}_q$ defined in Subsection \ref{geometry}.

\section{Proof of Theorem \ref{thm:spspr}} \label{sec:SpecialSpreads1}

Let $Q^+(5,q)$ be the Klein quadric in $\PG(5,q)$, $q$ odd, and let $\zeta$ be the orthogonal polarity of $\PG(5,q)$ associated to $Q^+(5,q)$. Let $\beta$ be a solid of $\PG(5,q)$ intersecting $Q^+(5,q)$ is an elliptic quadric, denoted here by $Q^-(3,q)$. Let $\alpha$ be a plane of $\PG(5,q)$ contained in $\beta$ intersecting $Q^-(3,q)$ (and hence also $Q^+(5,q)$) in an irreducible conic, denoted here by $Q_1(2,q)$.

Denote by $\zeta'$ the orthogonal polarity of $\beta$ associated to $Q^-(3,q)$. As $\alpha$ is a nontangent plane with respect to $Q^-(3,q)$, $z := \alpha^{\zeta'}$ is a point of $\beta$ not belonging to $Q^-(3,q)$ and hence also not to $Q^+(5,q)$. That means that $z^\zeta$ is a nontangent hyperplane with respect to $Q^+(5,q)$, thus intersecting $Q^+(5,q)$ is a parabolic quadric, denoted here by $Q(4,q)$.

As $z = \alpha^{\zeta'}$ all lines through $z$ containing a point of the irreducible conic $Q_1(2,q) = \alpha \cap Q^+(5,q)$ are tangent lines with respect to $Q^+(5,q)$, implying that $\alpha \subseteq z^\zeta$, or equivalently that $z \in \alpha^\zeta$. As $\alpha \cap Q^+(5,q)$ is an irreducible conic of $\alpha$ (namely $Q_1(2,q)$), also $\alpha^\zeta \cap Q^+(5,q)$ is an irreducible conic of $\alpha^\zeta$, denoted here by $Q_2(2,q)$. The point $z$ of the carrying plane $\alpha^\zeta$ of $Q_2(2,q)$ does not belong to $Q_2(2,q)$ and so must be either an interior point or an exterior point with respect to $Q_2(2,q)$. We will show that $z$ is an interior point with respect to $Q_2(2,q)$.

As $\alpha^\zeta$ is disjoint from $\alpha \subseteq z^\zeta$, it cannot be contained in $z^\zeta$ and so intersects $z^\zeta$ in a line $L$. We prove that $L$ is disjoint from $Q(4,q)$. Suppose to the contrary that $u \in L \cap Q(4,q)$. As $u \in z^\zeta$, we have $z \in u^\zeta$ and as $u \in \alpha^\zeta$, we also have $\alpha \subseteq u^\zeta$. So, $u^\zeta$ contains $\langle z,\alpha \rangle=\beta$. Now, $u^\zeta$ is a tangent hyperplane with respect to $Q^+(5,q)$ and thus intersects $Q^+(5,q)$ in a cone with vertex $u$ and as base a hyperbolic quadric in a hyperplane of $u^\zeta$ not containing $u$. But the same hyperplane would contain a solid (namely $\beta$) that intersects $Q^+(5,q)$ in an elliptic quadric (namely $Q^-(3,q)$), a clear contradiction. So, $L$ is indeed disjoint from $Q(4,q)$ and hence also from $Q^+(5,q)$.

Suppose now that $z$ is not an interior point with respect to $Q_2(2,q)$. Then $z$ is an exterior point with respect to $Q_2(2,q)$ and so there exists a tangent line $K$ in $\alpha^\zeta$ through $z$ intersecting $Q_2(2,q)$ in a point $v$. As $K$ is a tangent line, the point $v$ must belong to $z^\zeta$ and hence also to $z^\zeta \cap \alpha^\zeta = L$, in contradiction with $L \cap Q^+(5,q) = \emptyset$.

As $z$ is an interior point with respect to $Q_2(2,q)$, there are $\frac{q+1}{2}$ lines in $\alpha^\zeta$ through $z$ intersecting $Q_2(2,q)$ and hence $Q^+(5,q)$ in exactly two points. In $\beta$ itself, there are $q+1$ lines through $z$ intersecting $Q^-(3,q)$ in one point (namely the $q+1$ lines through $z$ and a point of $Q_1(2,q)$) and hence $\frac{|Q^-(3,q)| -(q+1)}{2}=\frac{q^2-q}{2}$ lines through $z$ intersecting $Q^-(3,q)$ (and hence also $Q^+(5,q)$) in two points. As $\alpha^\zeta \cap \alpha = \emptyset$, $\alpha^\zeta \cap \beta = \{ z \}$ and so we obtain $\frac{q^2+1}{2} = \frac{q+1}{2} + \frac{q^2-q}{2}$ lines through $z$ intersecting $Q^+(5,q)$ in exactly two points. The $q^2+1 = \frac{q^2+1}{2} \cdot 2$ points of $Q^+(5,q)$ we obtain in this way are exactly the points of $(Q^-(3,q) \setminus Q_1(2,q)) \cup Q_2(2,q)$. Note that $Q_1(2,q) = \alpha \cap Q^+(5,q)$ and $Q_2(2,q) = \alpha^\zeta \cap Q^+(5,q)$. As $Q^-(3,q)$ is an ovoid of $Q^+(5,q)$ also $O := (Q^-(3,q) \setminus Q_1(2,q)) \cup Q_2(2,q)$ is an ovoid of $Q^+(5,q)$, which is moreover disjoint from $Q(4,q)$. 

Now, let $\kappa$ be the Klein correspondence between the set of lines of $\PG(3,q)$ and the set of points of $Q^+(5,q)$. Then $\kappa^{-1}(Q(4,q))$ consists of all lines that are totally isotropic with respect to some symplectic polarity $\tau$ of $\PG(3,q)$. Denote by $W(q)$ the symplectic generalized quadrangle associated to $\tau$. The hyperbolic lines of $W(q)$ are the elements of $\tau^{-1}(Q^+(5,q) \setminus Q(4,q))$. The set $\kappa^{-1}(O)$ is thus a set of hyperbolic lines of $W(q)$. As $O$ is an ovoid of $Q^+(5,q)$, $\kappa^{-1}(O)$ is a line spread of $\PG(3,q)$.

We still need to show that for every $L \in \kappa^{-1}(O)$, we also have that $L^\tau$ belongs to $\kappa^{-1}(O)$. Put $L = \kappa^{-1}(y)$ for some $y \in O$. Then the line $zy$ intersects $O$ in a second point $y'$. As $Q(4,q) \subseteq z^\zeta$, we have that $y^\zeta \cap Q(4,q) = (y')^\zeta \cap Q(4,q) = yz^\zeta \cap Q(4,q)$, implying that the $(q+1)^2$ lines of $W(q)$ meeting $L = \kappa^{-1}(y)$ are exactly the $(q+1)^2$ lines of $W(q)$ meeting $L' = \kappa^{-1}(y')$, i.e. $\kappa^{-1}(y') = \kappa^{-1}(y)^\tau = L^\tau$ must indeed belong to $\kappa^{-1}(O)$.
 
The set $\kappa^{-1}(O)$ thus consists of hyperbolic lines of $W(q)$ and can be partitioned in pairs of the form $\{ L,L^\tau \}$.

\section{Proof of Theorem \ref{thm:2spspr}} \label{sec:SpecialSpreads2}

\noindent Let $W(q)$ with $q$ odd be the symplectic generalized quadrangle whose collinearity graph is $Sp(4,q)$. Let $S_1$ and $S_2$ be two special spreads of $\PG(3,q)$. We denote the complement of $\Gamma_{q,S_i}$ by $\overline{\Gamma_{q,S_i}}$.

\begin{lemma} \label{lem1bis}
For every $i \in \{ 1,2 \}$ and every point $x$ of $\PG(3,q)$, the local graph of $\overline{\Gamma_{q,S_i}}$ in the vertex $x$ is the disjoint union of $q+1$ cliques of size $q-1$.
\end{lemma}
\begin{proof}
Let $l$ be the unique element of $S_i$ containing $x$, and let $l^\perp$ be the hyperbolic line orthogonal to $l$. In $Sp(4,q)$, the vertices adjacent to $x$ are the vertices of the set $x^\zeta \setminus \{ x \}$, where $\zeta$ is the symplectic polarity of $\PG(3,q)$ associated to $W(q)$. The induced subgraph of $Sp(4,q)$ on the set $x^\zeta \setminus \{ x \}$ is the disjoint union of $q+1$ cliques of size $q$, and each such clique has the form $k \setminus \{ x \}$ where $k$ is a line of $x^\zeta$ through $x$. Now, each such line $k$ intersects $l^\perp$ in a unique point. In $\overline{\Gamma_{q,S_i}}$, the vertex $x$ is no longer adjacent with the vertices of $l^\perp$ and so the local graph of $\overline{\Gamma_{q,S_i}}$ in the vertex $x$ must be the disjoint union of $q+1$ cliques of size $q-1$.
\end{proof}

\medskip \noindent For every vertex $x$ of $\overline{\Gamma_{q,S_i}}$ and every clique $C$ of the local graph of $\overline{\Gamma_{q,S_i}}$ in the vertex $x$, the set $\{ x \} \cup C$ is called a {\em truncated symplectic line} of $\overline{\Gamma_{q,S_i}}$. From the proof of Lemma \ref{lem1bis}, it is clear that every truncated symplectic line of $\overline{\Gamma_{q,S_i}}$ is of the form $k \setminus \{ y \}$, where $k$ is a line of $W(q)$ and $y \in k$. In fact, the following can be proved.

\begin{lemma} \label{lem2bis}
The truncated symplectic lines of $\overline{\Gamma_{q,S_i}}$ are precisely the sets of the form $k \setminus \{ y \}$ where $k$ is a line of $W(q)$ and $y$ is some point of $k$.
\end{lemma}
\begin{proof}
It remains to show that every line of this form is a truncated symplectic line of $\overline{\Gamma_{q,S_i}}$. Let $l$ be the unique element of $S_i$ containing $y$. Then the line $k$ of $W(q)$ intersects the line $l^\perp$ in a unique point $x$. The set $k \setminus \{ x,y \}$ is then a maximal clique in the local graph of $\overline{\Gamma_{q,S_i}}$ in the vertex $x$. The associated truncated line for this clique is $(k \setminus \{ x,y \}) \cup \{ x \} = k \setminus \{ y \}$.
\end{proof}

\medskip \noindent The following is a consequence of Lemma \ref{lem2bis}.

\begin{lemma} \label{lem3}
The relation $R_i$, $i \in \{ 1,2 \}$, on the set of truncated lines of $\overline{\Gamma_{q,S_i}}$ defined by $(t_1,t_2) \in R_i$ if and only if either $t_1=t_2$ or $|t_1 \cap t_2|=q-1$ is an equivalence relation. There is a bijection between the set of equivalence classes and the lines of $W(q)$, namely if $k$ is a line of $W(q)$, then the set of subsets of size $q$ of $k$ is an equivalence class for the relation $R_i$.  
\end{lemma}

\medskip \noindent For an equivalence relation $C$ of the relation $R_i$, let $L_C$ denote the union of its elements and let $\mathcal{L}_i$ denote the set of all $L_C$'s, where $C$ is an equivalence class of the relation $R_i$. The elements of $\mathcal{L}_i$ are called the {\em symplectic lines} of $\overline{\Gamma_{q,S_i}}$.
 
\medskip \noindent The following is an immediate consequence of Lemma \ref{lem3}.
 
\begin{corollary} \label{co4}
For every $i \in \{ 1,2 \}$, $\mathcal{L}_i$ is the set of lines of $W(q)$.
\end{corollary} 
 
\medskip \noindent For every $i \in \{ 1,2 \}$, let $\mathcal{G}_i$ be the point-line geometry with line set $\mathcal{L}_i$ defined on the point set of $\PG(3,q)$. By Corollary \ref{co4}, $\mathcal{G}_i \cong W(q)$. Let $\mathcal{P}_i$ denote the set of all pairs $(u,v)$ of distinct points of $\PG(3,q)$ such that $u$ and $v$ are adjacent in $\mathcal{G}_i$, but not in $\overline{\Gamma_{q,S_i}}$. For every point $u$ of $\PG(3,q)$ and every $i \in \{ 1,2 \}$, let $T^{(i)}_u$ denote the set of all vertices $v$ such that $(u,v) \in \mathcal{P}_i$. Put $\mathcal{T}_i := \{ T^{(i)}_u \, | \, u \mbox{ is a point of } \PG(3,q) \}$. We call $\mathcal{T}_i$ the set of {\em hyperbolic lines} of $\overline{\Gamma_{q,S_i}}$.
 
\begin{lemma} \label{lem5bis}
Let $i \in \{ 1,2 \}$.
\begin{enumerate}
\item[$(1)$] The set $\mathcal{P}_i$ consists of all pairs $(u,v)$ where the respective lines of $S_i$ containing $u$ and $v$ are orthogonal.
\item[$(2)$] The hyperbolic lines of $\overline{\Gamma_{q,S_i}}$ are precisely the elements of $S_i$.
\end{enumerate}
\end{lemma}
\begin{proof}
Claim (1) is implied by the definitions of $\overline{\Gamma_{q,S_i}}$ and $\mathcal{P}_i$, along with the fact that $\mathcal{G}_i \cong W(q)$. Claim (2) is implied by (1).
\end{proof}

\medskip \noindent Note that the symplectic and hyperbolic lines of each $\overline{\Gamma_{q,S_i}}$, $i 
\in \{ 1,2 \}$, have been defined solely in terms of the adjacency relation of $\overline{\Gamma_{q,S_i}}$. This fact will be implicitly used in the proof of the following proposition.

\begin{proposition} \label{prop6}
The isomorphisms between the graphs $\overline{\Gamma_{q,S_1}}$ and $\overline{\Gamma_{q,S_2}}$ are precisely those automorphisms of $W(q)$ that map $S_1$ to $S_2$. 
\end{proposition}
 \begin{proof}
 From the definitions of the graphs $\overline{\Gamma_{q,S_1}}$ and $\overline{\Gamma_{q,S_2}}$, it immediately follows that every automorphism of $W(q)$ that maps $S_1$ to $S_2$ is also an isomorphism between $\overline{\Gamma_{q,S_1}}$ and $\overline{\Gamma_{q,S_2}}$. Conversely, suppose that $\theta$ is an isomorphism between $\overline{\Gamma_{q,S_1}}$ and $\overline{\Gamma_{q,S_2}}$. Then $\theta$ is a permutation of the point set of $\PG(3,q)$, i.e. of the point set of $W(q)$. The map $\theta$ must map the symplectic lines of $\overline{\Gamma_{q,S_1}}$ to the symplectic lines of $\overline{\Gamma_{q,S_2}}$. By Corollary \ref{co4}, we then know that $\theta$ is an automorphism of $W(q)$. As $\theta$ also maps the hyperbolic lines of $\overline{\Gamma_{q,S_1}}$ to the hyperbolic lines of $\overline{\Gamma_{q,S_2}}$, we know by Lemma \ref{lem5bis}(2) that $\theta$ maps $S_1$ to $S_2$.
\end{proof}

\medskip \noindent The following is a consequence of Proposition \ref{prop6}.

\begin{corollary}
\begin{enumerate}
\item[$(1)$] The automorphism group of $\Gamma_{q,S_i}$, $i \in \{ 1,2 \}$, consists of those automorphisms of $W(q)$ that fix the spread $S_i$.   
\item[$(2)$] The graphs $\Gamma_{q,S_1}$ and $\Gamma_{q,S_2}$ are isomorphic if and only if the spreads $S_1$ and $S_2$ of $W(q)$ are equivalent. 
\end{enumerate}
\end{corollary}

\section{The special spreads in $\PG(3,q)$, $q \in \{ 3,5,7 \}$} \label{sec:SpecialSpreads1.5}

Let $X = \{ 1,2,\ldots,n \}$ with $n \in \mathbb{N}$ be a set and $\mathcal{U}$ a set of subsets of $X$ whose union equals $X$. With the aid of the SageMath command $\mathrm{DLXCPP(\mathcal{U})}$ we can then determine all partitions of $X$ using only subsets of $\mathcal{U}$. We have used this command to find (all) partitions of the point set of $W(q)$, $q \in \{ 3,5,7 \}$, in subsets on which the induced subgraphs are isomorphic to $K_{q+1,q+1}$. We implemented this in two ways.

\medskip \noindent \textbf{Implementation 1:} The automorphism group $G_q$ of $W(q)$, $q \in \{ 3,5,7 \}$, is a group of type $PSp(4,q):2$. A model for the action of $G_q$ on the the point set of $W(q)$ can be found with the GAP command 
\[ \verb|g:=AllPrimitiveGroups(DegreeOperation,(q+1)*(q^2+1))[N(q)]|, \] 
where \verb|N(q)| equals 2 for $q \in \{ 3,5,7 \}$. In this GAP model, the point set $X_q$ coincides with $\{ 1,2,\ldots,(q+1)(q^2+1) \}$. We can subsequently compute stabilizers of (pairs of) points and their orbits on the point set. Selecting and/or merging some of these orbits allowed us to determine the set $\mathcal{L}_q$ of all lines, the set $\mathcal{H}_q$ of all hyperbolic lines and the set $\mathcal{U}_q$ of all point sets on which the induced subgraphs are isomorphic to $K_{q+1,q+1}$. We thus need to partition $X_q$ in subsets from $\mathcal{U}_q$. This goal can be achieved with the SageMath command \verb|DLXCPP| after replacing the pair $(X_q,\mathcal{U}_q)$ with an equivalent pair $(X,\mathcal{U})$ where $X$ equals $\{ 0,1,\ldots,(q+1)(q^2+1)-1 \}$. We computationally succeeded in this goal for $q \in \{ 3,5,7 \}$. In order to reduce the computations, we however modified this procedure for $q=5$ and $q=7$. As $G_q$ acts transitively on the elements of $\mathcal{U}_q$, we assumed that a given element $U$ of $\mathcal{U}_q$ (chosen in advance) belonged to the partition for $q=5$. For $q=7$, we first observed that $G_q$ has three orbits on the set of pairs $(U_1,U_2)$, where $U_1$ and $U_2$ are two disjoint elements of $\mathcal{U}_q$. We then did three separate computations. For a representative $(U_1,U_2)$ of each of the three orbits, we computed all partitions containing $U_1$ and $U_2$. Having found all partitions (containing $U$ for $q=5$ and $U_1,U_2$ for $q=7$ and for each of the three representatives $(U_1,U_2)$), we translated these partitions to the original GAP model, and subsequently we checked the isomorphism between two partitions \verb|P1| and \verb|P2| using the GAP command
\[ \verb|RepresentativeAction(g,P1,P2,OnSetsSets) <> fail|. \]
This method was successful for $q=3$ and $q=5$, where we respectively found 1 and 2 examples, up to isomorphism. This command was not able to provide any answers for $q=7$. Verifying whether two partitions are isomorphic, with each partition consisting of 25 sets of size 16 seemed computationally too hard. We therefore also implemented an alternative model of our problem, where the partitions correspond to sets of elements rather than to sets of sets of elements. In this way, we were also able to deal with the case $q=7$. 

\medskip \noindent \textbf{Implementation 2:} The point-line dual of the generalized quadrangle $W(q)$ is the generalized quadrangle $Q(4,q)$. The lines of $W(q)$ thus correspond to the points of $Q(4,q)$ and the points of $W(q)$ correspond to the lines of $Q(4,q)$. The pairs $\{ K,K^\zeta \}$ of orthogonal hyperbolic lines of $W(q)$ correspond bijectively to the hyperbolic quadrics of type $Q^+(3,q)$ on $Q(4,q)$ (obtained by intersecting $Q(4,q)$ with a suitable nontangent hyperplane). Such a hyperbolic quadric has two partitions in lines, the so-called {\em reguli}. The lines of one of these reguli correspond to the points of $K$, while the lines of the other reguli correspond to the points of $K^\zeta$. As the point sets of $W(q)$ on which the induced subgraphs are isomorphic to $K_{q+1,q+1}$ are precisely the elements of $\mathcal{U}_q$, there is thus a natural bijective correspondence between the elements of $\mathcal{U}_q$ and the nontangent hyperplanes intersecting $Q(4,q)$ in hyperbolic quadrics.

Now, let $\zeta'$ denote the orthogonal polarity of $\mathrm{PG}(4,q)$ associated to $Q(4,q)$. A point $x \in \mathrm{PG}(4,q) \setminus Q(4,q)$ is called {\em hyperbolic} if $x^{\zeta'} \cap Q(4,q)$ is a hyperbolic quadric and is called {\em elliptic} if $x^{\zeta'} \cap Q(4,q)$ is an elliptic quadric. The elements $U \in \mathcal{U}_q$ thus bijectively correspond to the hyperbolic points $x_U$ of $Q(4,q)$; if $\alpha_U$ is the nontangent hyperplane corresponding to $U$ (in the above sense), then $x_U := \alpha_U^{\zeta'}$. Note that two elements $U_1,U_2 \in \mathcal{U}_q$ are disjoint if and only if the hyperbolic quadrics $\alpha_{U_1} \cap Q(4,q)$ and $\alpha_{U_2} \cap Q(4,q)$ have no lines in common. 

\begin{lemma} \label{lem1}
Let $U_1$ and $U_2$ be two distinct elements of $\mathcal{U}_q$. Then the following are equivalent:
\begin{enumerate}
\item[$(1)$] $U_1$ and $U_2$ are not disjoint;
\item[$(2)$] there is a line $L$ of $Q(4,q)$ such that the points $x_{U_1}$ and $x_{U_2}$ belong to the plane $L^{\zeta'}$; 
\item[$(3)$] $x_{U_1} x_{U_2}$ is a tangent line.
\end{enumerate}
If one of these conditions is satisfied, then $|U_1 \cap U_2|=2$.
\end{lemma}
\begin{proof}
Suppose $U_1$ and $U_2$ are not disjoint. Let $y$ be a common point of $U_1$ and $U_2$, and let $L$ be the line of $Q(4,q)$ corresponding to $y$. Then $L$ is common line of the $Q^+(3,q)$-quadrics $x_{U_1}^{\zeta'} \cap Q(4,q)$ and $x_{U_2}^{\zeta'} \cap Q(4,q)$, implying that both $x_{U_1}$ and $x_{U_2}$ are contained in the plane $L^{\zeta'}$. So, (1) implies (2).

If there exists a line $L$ of $Q(4,q)$ such that $x_{U_1}$ and $x_{U_2}$ are contained in $L^{\zeta'}$, then the fact that $L^{\zeta'} \cap Q(4,q) = L$ implies that $x_{U_1} x_{U_2}$ is a tangent line. So, (2) implies (3).

Suppose now that $x_{U_1} x_{U_2}$ is a tangent line. There are two possibilities for the intersection of the $Q^+(3,q)$-quadrics $x_{U_1}^{\zeta'} \cap Q(4,q)$ and $x_{U_2}^{\zeta'} \cap Q(4,q)$. Either it is an irreducible conic in a plane $\alpha$ or it is the union of two lines. In the former case, $\alpha^{\zeta'}$ would be a line containing $x_{P_1}$ and $x_{P_2}$ that is either external or secant with respect to $Q(4,q)$, a clear contradiction. So, the intersection is the union of two lines. But then $U_1$ and $U_2$ have two points in common, and so are not disjoint.
\end{proof}

\bigskip \noindent Now, consider the following point-line geometry $\mathcal{S}_q$. The point set of $\mathcal{S}_q$ coincides with the set $\mathcal{P}_q$ of hyperbolic points of $\mathrm{PG}(4,q)$, and the lines of $\mathcal{S}_q$ are the nonempty intersections of $\mathcal{P}_q$ with the planes of $\mathrm{PG}(4,q)$ that intersect $Q(4,q)$ in a line, with incidence being containment. The geometry $\mathcal{S}_q$ has $|\mathcal{P}_q|=\frac{1}{2}q^2(q^2+1)$ points and $(q+1)(q^2+1)$ lines. The following clearly holds.

\begin{lemma} \label{lem2}
If $\mathcal{U} \subseteq \mathcal{U}_q$, then $\mathcal{U}$ forms a partition of the point set of $W(q)$ if and only if the set $\{ x_U \, | \, U \in \mathcal{U} \}$ is an ovoid of $\mathcal{S}_q$ (being a set of points of $\mathcal{S}_q$ having a unique point in common with each line of $\mathcal{S}_q$). 
\end{lemma}
\begin{proof}
Via the bijective correspondence between the elements $U \in \mathcal{U}_q$ and the points $x_U \in \mathcal{P}_q$, each line of $\mathcal{S}_q$ corresponds to a set of all elements of $\mathcal{U}_q$ containing a given point of $W(q)$. The claim follows.
\end{proof}

\bigskip \noindent For every hyperbolic point $x$ of $Q(4,q)$, let $\mathcal{U}_x$ denote the set of all planes of $\mathrm{PG}(4,q)$ through $x$ that intersect $Q(4,q)$ in a line. Let $X_q'$ denote the set of all planes intersecting $Q(4,q)$ in a line, and put $\mathcal{U}_q' := \{ \mathcal{U}_x \, | \, x \in \mathcal{P}_q \}$. The following is then an immediate consequence of Lemma \ref{lem2}.

\begin{corollary} \label{co3}
If $\mathcal{U} \subseteq \mathcal{U}_q$, then $\mathcal{U}$ forms a partition of the point set of $W(q)$ if and only if the set $\{ \mathcal{U}_{x_U} \, | \, U \in \mathcal{U} \}$ is a partition of $X_q'$ (using only elements of $\mathcal{U}_q'$).
\end{corollary}

\medskip \noindent In view of this corollary, it thus suffices to find all partitions of $X_q'$ in subsets belonging to the set $\mathcal{U}_q'$. This can in principle be achieved with the aid of the SageMath command \verb|DLXCPP| (if not computationally too hard).

\medskip \noindent The above requires that we implement a computer model of the geometry $\mathcal{S}_q$. The points of $\mathcal{S}_q$ correspond to the elements of $\mathcal{U}_q$ and the lines of $\mathcal{S}_q$ corresponds to a set of elements of $\mathcal{U}_q$ containing a given point of $W(q)$. The automorphism group of $W(q)$ therefore acts in a natural way as a group $G_q'$ of automorphisms of $\mathcal{S}_q$. A model for the group action of $G_q'$ on the points of $\mathcal{S}_q$ can be implemented with the GAP command
\[ \verb|g:=AllPrimitiveGroups(DegreeOperation,1/2*q^2*(q^2+1))[M(q)]|, \] 
where \verb|M(q)| equals $5$ if $q=3$, $2$ if $q=5$ and $19$ if $q=7$. Again by computing stabilizers of (pairs of) points and their orbits, it is possible to reconstruct all lines of $\mathcal{S}_q$ by selecting and/or combining some of these orbits. By means of Corollary \ref{co3} and the SageMath command \verb|DLXCPP| we can then find all partitions of $X_q'$ in elements of $\mathcal{U}_q'$. Each such partition $P$ thus corresponds to a subset $\{ x \in \mathcal{P}_q \, | \, \mathcal{U}_x \in P \}$ of $\mathcal{P}_q$. If we have two such partitions \verb|P1| and \verb|P2| with corresponding subsets \verb|S1| and \verb|S2| of $\mathcal{P}_q$, then we can verify whether they are equivalent under the group $G_q'$ with the aid of the following command:
\[ \verb|RepresentativeAction(g,S1,S2,OnSets) <> fail|. \]
This procedure easily worked for $q$ equal to 3, 5 and 7. This lead to the following conclusions:
\begin{itemize}
\item For $q=3$, there is up to $G_3'$-equivalence a unique set $S \subseteq \mathcal{P}_3$ for which $\{ \mathcal{U}_x \, | \, x \in S \}$ is a partition of $X_3'$. In Table \ref{tab1}, we have mentioned some information on its stabilizer (inside $G_3'$), such as the size and the structure. 
\item For $q=5$, there are up to $G_5'$-equivalence two sets $S \subseteq \mathcal{P}_5$ for which $\{ \mathcal{U}_x \, | \, x \in S \}$ is a partition of $X_5'$. In Table \ref{tab2}, we have mentioned some information on their stabilizers (inside $G_5'$).
\item For $q=7$, there are up to $G_7'$-equivalence two sets $S \subseteq \mathcal{P}_7$ for which $\{ \mathcal{U}_x \, | \, x \in S \}$ is a partition of $X_7'$. In Table \ref{tab3}, we have mentioned some information on their stabilizers (inside $G_7'$).
\end{itemize}
The following is a consequence of Corollary \ref{co3} and the above. 

\begin{theorem} \label{theo4}
\begin{enumerate}
\item[$(1)$] Up to isomorphism, there is a unique partition of the point set of $W(3)$ in subsets belonging to $\mathcal{U}_3$.
\item[$(2)$] Up to isomorphism, there are two partitions of the point set of $W(5)$ in subsets belonging to $\mathcal{U}_5$.
\item[$(3)$] Up to isomorphism, there are $14$ partitions of the point set of $W(7)$ in subsets belonging to $\mathcal{U}_7$.
\end{enumerate}
The structures of the stabilizers of these partitions (inside $G_q$ for $q \in \{ 3,5,7 \}$) along with their sizes can also be found in Tables $\ref{tab1}$, $\ref{tab2}$ and $\ref{tab3}$.
\end{theorem}

\begin{table}
\begin{center}
\begin{tabular}{|c|c|c|c|}
\hline
Example & Structure stabilizer & Size & Characteristic \\
\hline \hline
1 & $((C_2 \times C_2 \times C_2 \times C_2) : A_5) : C_2$ & 1920 & [ 0,10 ] \\
\hline
\end{tabular}
\end{center}
\caption{The partitions for $q=3$} \label{tab1}
\end{table}

\begin{table}
\begin{center}
\begin{tabular}{|c|c|c|c|}
\hline
Example & Structure stabilizer & Size & Characteristic \\
\hline \hline
1 & $((((SL(2,3) : C_2) : C_2) : C_3) : C_2) : C_2$ & 1152 & [ 0, 48, 30 ] \\
\hline
2 & $C_2 \times S_5 \times S_3$ & 1440 & [ 0, 33, 45 ] \\
\hline
\end{tabular}
\end{center}
\caption{The partitions for $q=5$} \label{tab2}
\end{table}

\begin{table}
\begin{center}
\begin{tabular}{|c|c|c|c|}
\hline
Example & Structure stabilizer & Size & Characteristic \\
\hline
\hline
1 & $((((C_2 \times D_8) : C_2) : C_3) : C_2) : C_2$ & 384 & [ 0, 156, 96, 48 ] \\
\hline
2 & $S_6$ & 720 & [ 0, 45, 90, 165 ] \\
\hline
3 & $C_2 \times S_4 \times S_3$ & 288 & [ 0, 57, 48, 195 ] \\
\hline
4 & $(D_{16} \times D_{16}) : C_2$ & 512 & [ 0, 72, 108, 120 ] \\
\hline
5 & $((C_2 \times C_2 \times C_2 \times C_2) : A_5) : C_2$ & 1920 & [ 0, 120, 140, 40 ] \\
\hline
6 & $(D_{16} \times D_{16}) : C_2$ & 512 & [ 0, 128, 116, 56 ] \\
\hline
7 & $D_{16} \times (PSL(3,2) : C_2)$ & 5376 & [ 0, 84, 128, 88  \\
\hline
8 & $C_2 \times C_2 \times S_4$ & 96 & [ 0, 141, 62, 97 ] \\
\hline
9 & $(D_8 \times D_8) : C_2$ & 128 & [ 0, 132, 72, 96 ] \\
\hline
10 & $((C_2 \times C_2 \times C_2 \times C_2) : C_5) : C_2$ & 160 & [ 0, 60, 100, 140 ] \\
\hline
11 & $(((((SL(2,3) : C_2) : C_2) : C_3) : C_2) : C_2) : C_2$ & 2304 & [ 0, 96, 132, 72 ] \\
\hline
12 & $(D_8 \times D_8) : C_2$ & 128 & [ 0, 128, 36, 136 ] \\
\hline
13 & $C_{25} : C_4$ & 100 & [ 0, 150, 50, 100 ] \\
\hline
14 & $C_2 \times S_5$ & 240 & [ 0, 120, 60, 120 ] \\
\hline
\end{tabular}
\end{center}
\caption{The partitions for $q=7$} \label{tab3}
\end{table}

\medskip \noindent We now wish to determine some structural properties of the partitions which uniquely determine their isomorphism classes. In this way, we can easily make identification between isomorphic partitions in different models of $W(q)$. Also, it will allow us to make such identifications in the same model where the isomorphism test with other methods will fail (for instance, because it is computationally too hard as in the case $q=7$ in Implementation 1, where the GAP command gave no definite answer, see above). We first need to do some preparatory work.   

\medskip \noindent Let $V$ be a 5-dimensional vector space over the field $\mathbb{F}_q$ for which $\mathrm{PG}(4,q)=\mathrm{PG}(V)$, let $\mathcal{Q}: V \to \mathbb{F}_q$ be a quadratic form on $V$ such that $Q(4,q)$ consists of all points $\langle \bar v \rangle$ of $\mathrm{PG}(4,q)$ for which $\mathcal{Q}(\bar v) = 0$ and let $\mathcal{B}: V \times V \to \mathbb{F}_q$ be the bilinear form on $V$ defined by $\mathcal{Q}$, i.e. $\mathcal{B}(\bar v_1,\bar v_2) = \mathcal{Q}(\bar v_1 + \bar v_2) - \mathcal{Q}(\bar v_1) - \mathcal{Q}(\bar v_2)$ for all $\bar v_1,\bar v_2 \in V$. 

\begin{lemma} \label{lem5}
Let $\bar v_1,\bar v_2$ be two linearly independent vectors of $V$ such that $\mathcal{Q}(\bar v_1) \not= 0 \not= \mathcal{Q}(\bar v_2)$ and $\mathcal{B}(\bar v_1,\bar v_2)=0$. Let $L$ be the line of $\mathrm{PG}(4,q)$ containing $\langle \bar v_1 \rangle$ and $\langle \bar v_2 \rangle$. Then $L$ is a secant line with respect to $Q(4,q)$ if and only if $-\frac{\mathcal{Q}(\bar v_1)}{\mathcal{Q}(\bar v_2)}$ is a square in $\mathbb{F}_q$ and an external line with respect to $Q(4,q)$ if and only if $-\frac{\mathcal{Q}(\bar v_1)}{\mathcal{Q}(\bar v_2)}$ is a nonsquare in $\mathbb{F}_q$.
\end{lemma}
\begin{proof}
Obviously, none of the points $\langle \bar v_1 \rangle$, $\langle \bar v_2 \rangle$ belongs to $Q(4,q)$. The points of $L \setminus \{ \langle \bar v_2 \rangle \}$ have the form $\langle \bar v_1 + \lambda \bar v_2 \rangle$, $\lambda \in \mathbb{F}_q$. The equation $\mathcal{Q}(\bar v_1 + \lambda \bar v_2) = \mathcal{Q}(\bar v_1) + \lambda^2 \mathcal{Q}(\bar v_2)=0$ has (necessarily two) solutions if and only if $-\frac{\mathcal{Q}(\bar v_1)}{\mathcal{Q}(\bar v_2)} $ is a square in $\mathbb{F}_q$. The claim follows.
\end{proof}

\medskip \noindent For every point $x \in \mathcal{P}_q$, we denote by $x^\perp$ the set of all $y \in \mathcal{P}_q$ for which $y \in x^{\zeta'}$. We can now choose the quadratic form $\mathcal{Q}$ such that $\mathcal{P}_q$ consists of all points $\langle \bar v \rangle$ of $\mathrm{PG}(4,q)$ such that $\mathcal{Q}(\bar v)$ is a nonzero square. A secant line contains exactly $\frac{q-1}{2}$ points of $\mathcal{P}_q$ and an external line contains exactly $\frac{q+1}{2}$ points of $\mathcal{P}_q$. Through each point of $\mathcal{P}_q$, there are $(q+1)^2$ tangent lines, $\frac{1}{2}\Big( (q+1)(q^2+1)-(q+1)^2 \Big) = \frac{1}{2}(q^3-q)$ secant lines and $(q^3+q^2+q+1)-\frac{1}{2}(q^3-q)-(q+1)^2=\frac{1}{2}(q^3-q)$ external lines. These facts in combination with Lemma \ref{lem5} allows us to draw the following conclusions.

\begin{corollary} \label{lem6}
For $q \equiv 1 \pmod{4}$, the set $\mathcal{P}_q \times \mathcal{P}_q$ can be partitioned into the following five subsets:
\begin{enumerate}
\item[$(a)$] the set of all pairs $(x,x)$ where $x \in \mathcal{P}_q$;
\item[$(b)$] the set of all pairs $(x,y) \in \mathcal{P}_q \times \mathcal{P}_q$ with $x \not= y$ and $xy$ is a tangent line (for any such $(x,y)$, we have $y \in x^{\zeta'}$);
\item[$(c)$] the set of all pairs $(x,y) \in \mathcal{P}_q \times \mathcal{P}_q$ with $x \not= y$, $y \in x^{\zeta'}$ and $xy$ is a secant line;
\item[$(d)$] the set of all pairs $(x,y) \in \mathcal{P}_q \times \mathcal{P}_q$ with $x \not= y$, $y \not\in x^{\zeta'}$ and $xy$ is a secant line;
\item[$(e)$] the set of all pairs $(x,y) \in \mathcal{P}_q \times \mathcal{P}_q$ with $x \not= y$ and $xy$ is an external line (for any such $(x,y)$, we have $y \not\in x^{\zeta'}$). 
\end{enumerate}
The respective sizes of these subset of $\mathcal{P}_q \times \mathcal{P}_q$ are $|\mathcal{P}_q|$, $|\mathcal{P}_q| \cdot (q-1)(q+1)^2$, $|\mathcal{P}_q| \cdot \frac{1}{2}(q^3-q)$, $|\mathcal{P}_q| \cdot \frac{1}{2}(q^3-q)\frac{q-5}{2}$ and $|\mathcal{P}_q| \cdot \frac{1}{2}(q^3-q) \frac{q-1}{2}$. These numbers are mutually distinct. The set of all $(x,y) \in \mathcal{P}_q \times \mathcal{P}_q$ for which $y \in x^\perp$ can be obtained as the union of two of these subsets, namely the ones in $(b)$ and $(c)$.
\end{corollary}

\medskip \noindent Suppose $q \equiv 1 \pmod{4}$. Having available the partition of $\mathcal{P}_q \times \mathcal{P}_q$ in subsets as in Corollary \ref{lem6}, we can determine the nature of these subsets, solely based on their sizes, in particular, based on the information of the sizes alone, we can determine the set $x^\perp$ for every point $x \in \mathcal{P}_q$. This method of determining these subsets of $\mathcal{P}_q \times \mathcal{P}_q$ and the sets $x^\perp$ for points $x \in \mathcal{P}_q$ can be implemented in our GAP computer model of the geometry $\mathcal{S}_5$. 

\begin{corollary} \label{lem7}
For $q \equiv 3 \pmod{4}$, the set $\mathcal{P}_q \times \mathcal{P}_q$ can be partitioned into the following five subsets:
\begin{enumerate}
\item[$(a)$] the set of all pairs $(x,x)$ where $x \in \mathcal{P}_q$;
\item[$(b)$] the set of all pairs $(x,y) \in \mathcal{P}_q \times \mathcal{P}_q$ with $x \not= y$ and $xy$ is a tangent line (for any such $(x,y)$, we have $y \in x^{\zeta'}$);
\item[$(c)$] the set of all pairs $(x,y) \in \mathcal{P}_q \times \mathcal{P}_q$ with $x \not= y$ and $xy$ is a secant line (for any such $(x,y)$, we have $y \not\in x^{\zeta'}$);
\item[$(d)$] the set of all pairs $(x,y) \in \mathcal{P}_q \times \mathcal{P}_q$ with $x \not= y$, $y \in x^{\zeta'}$ and $xy$ is an external line;
\item[$(e)$] the set of all pairs $(x,y) \in \mathcal{P}_q \times \mathcal{P}_q$ with $x \not= y$, $y \not\in x^{\zeta'}$ and $xy$ is an external line. 
\end{enumerate}
The respective sizes of these subset of $\mathcal{P}_q \times \mathcal{P}_q$ are $|\mathcal{P}_q|$, $|\mathcal{P}_q| \cdot (q-1)(q+1)^2$, $|\mathcal{P}_q| \cdot \frac{1}{2}(q^3-q) \frac{q-3}{2}$, $|\mathcal{P}_q| \cdot \frac{1}{2}(q^3-q)$ and $|\mathcal{P}_q| \cdot \frac{1}{2}(q^3-q) \frac{q-3}{2}$. With exception of the size $\frac{1}{2}(q^3-q) \frac{q-3}{2}$, which occurs twice, all these numbers are mutually distinct. The set of all $(x,y) \in \mathcal{P}_q \times \mathcal{P}_q$ for which $y \in x^\perp$ can be obtained as the union of two of these subsets, namely the ones in $(b)$ and $(d)$.
\end{corollary}

\medskip \noindent Suppose $q \equiv 3 \pmod{4}$. Having available the partition of $\mathcal{P}_q \times \mathcal{P}_q$ in subsets as in Corollary \ref{lem7}, we can determine the nature of some of these subsets, solely based on their sizes, namely the subsets mentioned in (a), (b) and (d), in particular, based on the information of the sizes alone, we can determine the set $x^\perp$ for every point $x \in \mathcal{P}_q$ (just as in the case where $q \equiv 1 \pmod{4}$). This method can be implemented in our GAP computer model of the geometry $\mathcal{S}_q$, $q \in \{ 3,7 \}$. We would now also like to find a method which allows us to determine all subsets (thus also those mentioned in (c) and (e)). For this, it suffices to give a criterion (again implementable in our GAP computer models of $\mathcal{S}_q$, $q \in \{ 3,7 \}$) to determine for any two distinct points $x$ and $y$ of $\mathcal{P}_q$ for which $xy$ is not a tangent line whether the line $xy$ is a secant or an external line. This will be achieved in the following lemma.

\begin{lemma} \label{lem8}
Let $x$ and $y$ be two distinct points of $\mathcal{P}_q$ such that $xy$ is not a tangent line. Then $xy$ is a secant line if and only if $|x^\perp \cap y^\perp| = \frac{q(q+1)}{2}$ and $xy$ is an external line if and only if $|x^\perp \cap y^\perp| = \frac{q(q-1)}{2}$. 
\end{lemma}
\begin{proof}
As the line $xy$ is not a tangent line, it is either a secant line or an external line, and so the plane $\alpha := xy^{\zeta'}$ intersects $Q(4,q)$ in an irreducible conic $C$. Let $E$ (respectively, $I$) denote the set of all points of $\alpha$ that are {\em exterior} (respectively, {\em interior}) with respect to $C$. Then $|E| = \frac{q(q+1)}{2}$ and $|I| = \frac{q(q-1)}{2}$. The set $x^\perp \cap y^\perp = \mathcal{P}_q \cap \alpha$ coincides with either $E$ or $I$ (and so contains either $\frac{q(q+1)}{2}$ or $\frac{q(q-1)}{2}$ points), see e.g. Lemma 2.7(1) of \cite{DB-Pa-Pr-Sa}. Let $z$ denote an arbitrary point of $E$ and denote by $u_1$ and $u_2$ the two points of $C$ such that $zu_1$ and $zu_2$ are tangent lines. Now, $Q'' := z^{\zeta'} \cap Q(4,q)$ is either a hyperbolic quadric (if $z \in \mathcal{P}_q$) or an elliptic quadric (if $z \not\in \mathcal{P}_q$) of $z^{\zeta'}$ and we denote by $\zeta''$ the orthogonal polarity of $z^{\zeta'}$ associated with $Q''$. Note that $u_1,u_2,x,y \in z^{\zeta'}$ and $u_1^{\zeta''} \cap u_2^{\zeta''} = xy$. If $Q''$ is a hyperbolic quadric, then $u_1^{\zeta''} \cap u_2^{\zeta''} \cap Q'' = xy \cap Q'' = xy \cap Q$ must be a set of two points, and if $Q''$ is an elliptic quadric, then $xy \cap Q$ must be empty. So, we have 
\[ |x^\perp \cap y^\perp| = \frac{q(q+1)}{2} \Leftrightarrow z \in \mathcal{P}_q  \Leftrightarrow Q'' \mbox{ is a hyperbolic quadric } \Leftrightarrow xy \mbox{ is a secant line}. \]
\end{proof}

\medskip \noindent We now take a closer look to the cases  where $q \in \{ 3,5,7 \}$.

For $q=3$, we know by Corollary \ref{lem7} that the set $\mathcal{P}_q \times \mathcal{P}_q$ can be partitioned in $\{ R_0,R_1,R_2 \}$, where $R_0$ consists of all pairs $(x,x)$ with $x \in \mathcal{P}_q$, $R_1$ consists of all pairs $(x,y) \in \mathcal{P}_q \times \mathcal{P}_q$ for which $x \not= y$ and $xy$ is a tangent line, and $R_2$ consists of all pairs $(x,y) \in \mathcal{P}_q \times \mathcal{P}_q$ for which $x \not= y$, $y \in x^{\zeta'}$ and $xy$ is an external line.

For $q=5$, we know by Corollary \ref{lem6} that the set $\mathcal{P}_q \times \mathcal{P}_q$ can be partitioned in $\{ R_0,R_1,R_2,R_3 \}$, where $R_0$ consists of all pairs $(x,x)$ with $x \in \mathcal{P}_q$, $R_1$ consists of all pairs $(x,y) \in \mathcal{P}_q \times \mathcal{P}_q$ for which $x \not= y$ and $xy$ is a tangent line, $R_2$ consists of all pairs $(x,y) \in \mathcal{P}_q \times \mathcal{P}_q$ for which $x \not= y$, $y \in x^{\zeta'}$ and $xy$ is a secant line, and $R_3$ consists of all pairs $(x,y) \in \mathcal{P}_q \times \mathcal{P}_q$ for which $x \not= y$, $y \not\in x^{\zeta'}$ and $xy$ is an external line.

For $q=7$, we know by Corollary \ref{lem7} that the set $\mathcal{P}_q \times \mathcal{P}_q$ can be partitioned in $\{ R_0,R_1,R_2,R_3,R_4 \}$, where $R_0$ consists of all pairs $(x,x)$ with $x \in \mathcal{P}_q$, $R_1$ consists of all pairs $(x,y) \in \mathcal{P}_q \times \mathcal{P}_q$ for which $x \not= y$ and $xy$ is a tangent line, $R_2$ consists of all pairs $(x,y) \in \mathcal{P}_q \times \mathcal{P}_q$ for which $x \not= y$, $y \not\in x^{\zeta'}$ and $xy$ is a secant line, $R_3$ consists of all pairs $(x,y) \in \mathcal{P}_q \times \mathcal{P}_q$ for which $x \not= y$, $y \in x^{\zeta'}$ and $xy$ is an external line, and $R_4$ consists of all pairs $(x,y) \in \mathcal{P}_q \times \mathcal{P}_q$ for which $x \not= y$, $y \not\in x^{\zeta'}$ and $xy$ is an external line.

Now put $i^\ast := 2$ if $q=3$, $i^\ast := 3$ if $q=5$ and $i^\ast := 4$ if $q=7$. For a subset $X$ of $\mathcal{P}_q$ with $q \in \{ 3,5,7 \}$, we define the {\em characteristic} of $X$ as the $i^\ast$-tuple $[N_1,N_2,\ldots,N_{i^\ast}]$, where $N_i$ with $i \in \{ 1,2,\ldots,i^\ast \} $ denotes the cardinality of $R_i \cap (X \times X)$. 

We have verified in a computational way that for $q \in \{ 3,5,7 \}$, there are $i^\ast$ orbits of $G_q'$ on the set of ordered pairs $(x,y)$, where $x$ and $y$ are two distinct elements of $\mathcal{P}_q$. This in combination with the preceding discussion implies the following.

\begin{corollary} \label{co9}
Let $q \in \{ 3,5,7 \}$. Then the $i^\ast$ orbits of $G_q'$ on the set of ordered pairs $(x,y) \in \mathcal{P}_q \times \mathcal{P}_q$ with $x \not= y$ are exactly the sets $R_1,R_2,\ldots,R_{i^\ast}$.
\end{corollary}

\medskip \noindent Using GAP, we have computed the $i^\ast$ orbits of $G_q'$ in the set of $\mathcal{P}_q \times \mathcal{P}_q$, for $q \in \{ 3,5,7 \}$. With the aid of Lemma \ref{lem8} and Corollaries \ref{lem6}, \ref{lem7} with their ensuing discussions, we can then explicitly determine which of these orbits correspond to $R_1,R_2,\ldots,R_{i^\ast}$. Once we have this information, we can determine the characteristic of any subset of $\mathcal{P}_q$. We have applied this to the representatives of the $G_q'$-equivalence classes of subsets $S \subseteq \mathcal{P}_q$, $q \in \{ 3,5,7 \}$, for which $\{ \mathcal{U}_x \, | \, x \in S \}$ is a partition of $X_q'$. This information has also been included in Tables \ref{tab1}, \ref{tab2} and \ref{tab3}. We this see that the $G_q'$-equivalence classes are uniquely determined by the characteristics of their representatives.

\medskip We can also use these characteristics to distinguish between the various isomorphism classes in our original model of the problem, where we need to partition the point set of $W(q)$ in elements of $\mathcal{U}_q$. This method might be effective if we do not have the group action available in our computer implementation of this model, or if the implemented group action is not able to provide any answers (as it was the case for $q=7$ in Implementation 1, see above). In order to determine the characteristics we need to solve the following problems for two distinct elements $U_1,U_2 \in \mathcal{U}_q$ with respective associated elements $x_{U_1},x_{U_2} \in \mathcal{P}_q$:
\begin{itemize}
\item[(a)] What conditions need to be satisfied by $U_1$ and $U_2$ for $x_{U_1}x_{U_2}$ to be a tangent, secant line or external line with respect to $Q(4,q)$?
\item[(b)] What conditions need to be satisfied by $U_1$ and $U_2$ for $x_{U_2}$ to belong to $x_{U_1}^{\zeta'}$?
\end{itemize}
Problem (a) will be solved in Lemmas \ref{lem10} and \ref{lem11} below. Problem (b) will be solved in Lemma \ref{lem13}. We first need to make a definition. For an element $U \in \mathcal{U}_q$, let $\mathcal{L}_U$ denote the set of $(q+1)^2$ symplectic lines intersecting $U$ in exactly two points, i.e. each of the two orthogonal hyperbolic lines of $U$ in exactly one point.

\begin{lemma} \label{lem10}
If $U_1$ and $U_2$ are two distinct elements of $\mathcal{U}_q$, then $|\mathcal{L}_{U_1} \cap \mathcal{L}_{U_2}| \in \{ q+1,2q+1 \}$. Moreover, $|\mathcal{L}_{U_1} \cap \mathcal{L}_{U_2}|=2q+1$ if and only if $x_{U_1}x_{U_2}$ is a tangent line, and if $|\mathcal{L}_{U_1} \cap \mathcal{L}_{U_2}|=q+1$, then the $q+1$ lines of $\mathcal{L}_{U_1} \cap \mathcal{L}_{U_2}$ are mutually disjoint.
\end{lemma}
\begin{proof}
A symplectic line $L$ belongs to $\mathcal{L}_{U_i}$, $i \in \{ 1,2 \}$, if and only if the point of $Q(4,q)$ corresponding to $L$ is contained in $X_{U_i}^{\zeta'}$. So, $|\mathcal{L}_{U_1} \cap \mathcal{L}_{U_2}|$ equals the size of the intersection of the two $Q^+(3,q)$-quadrics $x_{U_1}^{\zeta'} \cap Q(4,q)$ and $x_{U_2}^{\zeta'} \cap Q(4,q)$. The intersection is either an irreducible conic (containing $q+1$ points) or the union of two lines (containing $2q+1$ points). Hence, $|\mathcal{L}_{U_1} \cap \mathcal{L}_{U_2}| \in \{ q+1,2q+1 \}$. In the proof of Lemma \ref{lem1}, we have already seen that $x_{U_1}x_{U_2}$ is a tangent line if and only if $x_{U_1}^{\zeta'} \cap x_{U_2}^{\zeta'} \cap Q(4,q)$ is the union of two lines.

If $|\mathcal{L}_{U_1} \cap \mathcal{L}_{U_2}|=q+1$, then the $q+1$ points corresponding to the lines of $\mathcal{L}_{U_1} \cap \mathcal{L}_{U_2}$ form an irreducible conic. These points are therefore mutually noncollinear on $Q(4,q)$, implying that the lines in $\mathcal{L}_{U_1} \cap \mathcal{L}_{U_2}$ are mutually disjoint.
\end{proof}

\medskip \noindent If $K_1$, $K_2$ and $K_3$ are mutually disjoint lines of $\mathrm{PG}(3,q)$, there are exactly $q+1$ lines in $\mathrm{PG}(3,q)$ meeting $K_1$, $K_2$ and $K_3$, see e.g. Lemma 15.1.1 of \cite{Hi}. These $q+1$ lines are moreover mutually disjoint.

\begin{lemma} \label{lem11}
Let $U_1$ and $U_2$ be two distinct elements of $\mathcal{U}_q$ such that $|\mathcal{L}_{U_1} \cap \mathcal{L}_{U_2}|=q+1$. Let $L_1,L_2,L_3,L_1',L_2',L_3' \in \mathcal{L}_{U_1} \cap \mathcal{L}_{U_2}$ such that $L_1,L_2,L_3$ are mutually distinct as well as $L_1',L_2',L_3'$. Let $\mathcal{L}$ (respectively, $\mathcal{L}'$) denote the set of $q+1$ lines of $\mathrm{PG}(3,q)$ meeting $L_1,L_2,L_3$ (respectively, $L_1',L_2',L_3'$). Then the following hold:
\begin{enumerate}
\item[$(1)$] the singular lines contained in $\mathcal{L}$ are exactly the singular lines contained in $\mathcal{L}'$; 
\item[$(2)$] the number of singular lines contained in $\mathcal{L}$ (or $\mathcal{L}'$) equals $0$ or $2$;
\item[$(3)$] the number of singular lines contained in $\mathcal{L}$ (or $\mathcal{L}'$) equals $2$ if and only if $x_{U_1}x_{U_2}$ is a secant line;
\item[$(4)$] the number of singular lines contained in $\mathcal{L}$ (or $\mathcal{L}'$) equals $0$ if and only if $x_{U_1}x_{U_2}$ is an external line.
\end{enumerate}
\end{lemma}
\begin{proof}
As $|\mathcal{L}_{U_1} \cap \mathcal{L}_{U_2}|=q+1$, we know from Lemma \ref{lem10} that $x_{U_1}x_{U_2}$ is not a tangent line and so is either a secant line or an external line. In order to prove (1), (2), (3), (4), it suffices to prove that the singular lines (i.e. lines of $W(q)$) contained in $\mathcal{L}$ (or $\mathcal{L}'$) correspond to those points of $Q(4,q)$ that are contained on the line $x_{U_1}x_{U_2}$. 

Let $\alpha$ be the plane arising as intersection of the solids $x_{U_1}^{\zeta'}$ and $x_{U_2}^{\zeta'}$. In the proof of Lemma \ref{lem1}, we have seen that $\alpha \cap Q(4,q)$ must be an irreducible conic $C$. 

Now, let $x$ be an arbitrary point of $Q(4,q)$. The singular line of $W(q)$ corresponding to $x$ belongs to $\mathcal{L}$ if and only if $x$ is collinear on $Q(4,q)$ with the points $x_1$, $x_2$ and $x_3$ of $Q(4,q)$ corresponding to respectively $L_1$, $L_2$ and $L_3$ (these three points belong to $C$). This precisely happens when $x$ is collinear on $Q(4,q)$ with all points of $C$, i.e. if and only if $x \in \alpha^{\zeta'} \cap Q(4,q) = x_{U_1}x_{U_2} \cap Q(4,q)$. 

In a similar way, one shows that the singular line of $W(q)$ corresponding to $x$ belongs to $\mathcal{L}'$ if and only if $x \in  x_{U_1}x_{U_2} \cap Q(4,q)$. The claims follow.
\end{proof}

\begin{lemma} \label{lem12}
Let $L_1$ and $L_2$ be two disjoint singular lines of $W(q)$, and let $U \in \mathcal{U}_q$. If $x_1$ and $x_2$ are the points of $Q(4,q)$ corresponding to respectively $L_1$ and $L_2$, then $x_U \in x_1x_2$ if and only if each of the $q+1$ singular lines meeting $L_1$ and $L_2$ belongs to $\mathcal{L}_U$.
\end{lemma}
\begin{proof}
Note that as $L_1 \cap L_2 = \emptyset$, the points $x_1$ and $x_2$ are noncollinear on $Q(4,q)$. We have $x_U \in x_1x_2$ if and only if $x_1^{\zeta'} \cap x_2^{\zeta'} = (x_1x_2)^{\zeta'} \subset x_U^{\zeta'}$. The plane $(x_1x_2)^{\zeta'}$ intersects $Q(4,q)$ is a set of $q+1$ points  (an irreducible conic) corresponding to the $q+1$ singular lines of $W(q)$ meeting $L_1$ and $L_2$. The solid $x_U^{\zeta'}$ intersects $Q(4,q)$ in a set of $(q+1)^2$ points (a hyperbolic quadric) corresponding to the $(q+1)^2$ lines of $\mathcal{L}_U$. The claim follows.
\end{proof}

\begin{lemma} \label{lem13}
Let  $U_1$ and $U_2$ be two distinct elements of $\mathcal{U}$. Then $x_{U_2} \in x_{U_1}^{\zeta'}$ if and only if there exist two disjoint singular lines $K$ and $L$ in $\mathcal{L}_{U_1}$ such that all $q+1$ singular lines meeting $K$ and $L$ are contained in $\mathcal{L}_{U_2}$.   
\end{lemma}
\begin{proof}
The solid $x_{U_1}^{\zeta'}$ intersects $Q(4,q)$ in a hyperbolic quadric $Q$, whose $(q+1)^2$ points correspond to the $(q+1)^2$ lines contained in $\mathcal{L}_{U_1}$. We have that $x_{U_2} \in x_{U_1}^{\zeta'}$ if and only if there are two points in $x_{U_1}^{\zeta'} \cap Q(4,q)$ which are contained together with $x_{U_2}$ on the same line. The claim then follows from Lemma \ref{lem12}. 
\end{proof}

\medskip \noindent Based on Lemmas \ref{lem10}, \ref{lem11} and \ref{lem13}, we can now find for every subset $\mathcal{U}$ of $\mathcal{U}_q$, $q \in \{ 3,5,7 \}$, the characteristic of the corresponding set $\{ x_U \, | \, U \in \mathcal{U} \} \subseteq \mathcal{P}_q$ (allowing us to independently verify the information in the last column of Tables \ref{tab1}, \ref{tab2} and \ref{tab3} for the partitions of the point set of $W(q)$ in elements of $\mathcal{U}_q$). This also provides a method to verify whether two partitions \verb|P1| and \verb|P2| of the point set of $W(q)$ in elements of $\mathcal{U}_q$ are isomorphic. This verification requires a relatively small of amount of computations, and immediately provides answers for each of the prime powers $3$, $5$ and $7$, while in the previous method (Implementation 1) the command  \verb|RepresentativeAction(g,P1,P2,OnSetsSets) <> fail| was not even able to provide any answer for $q=7$.

\section{Further discussion}\label{sec:Discussion}

In this section we discuss the main results of this paper and formulate several open problems.

We have verified computationally that, for $q = 3,5$ and $7$, the projective space $\operatorname{PG}(3,q)$ has exactly $1,2$ and $14$ pairwise non-equivalent special spreads. We then verified that only one graph can be obtained by applying Theorem \ref{thm:DDG2} to the special spread for $q = 3$. Moreover, 12 and 16 graphs, respectively can be obtained by applying Theorem \ref{thm:DDG2} to the two special spreads for $q = 5$. Among these 28 graphs there are two pairs of isomorphic graphs. Thus, for $q = 5$, Theorem \ref{thm:DDG2} gives 26 pairwise non-isomorphic graphs in total. In case $q = 7$, we took one of the fourteen special spreads, got at least 6000 pairwise non-isomorphic graphs from Theorem \ref{thm:DDG2} and stopped the search.

We thus formulate the following open problems on special spreads.

\begin{problem}
Given an odd prime power $q$, how many pairwise non-equivalent special spreads does there exist in $\operatorname{PG}(3,q)$? In other words, given an odd prime power $q$, how many pairwise non-isomorphic graphs does Theorem \ref{thm:DDG1} produce?    
\end{problem}

\begin{problem}
Given an odd prime power $q$, how many pairwise non-isomorphic graphs does Theorem \ref{thm:DDG2} produce?    
\end{problem}

In a similar way, we formulate the following open problems on symplectic spreads.

\begin{problem}
Given an odd prime power $q$, how many pairwise non-equivalent symplectic spreads does there exist in $\operatorname{PG}(3,q)$?    
\end{problem}

\begin{problem}
Given an odd prime power $q$, how many pairwise non-isomorphic graphs does Theorem \ref{thm:DDG3} produce?    
\end{problem}

\begin{problem}
Given an odd prime power $q$, how many pairwise non-isomorphic graphs does Theorem \ref{thm:DDG4} produce?    
\end{problem}




Finally, we focus on equitable partitions of $Sp(4,q)$ and, in particular, equitable 2-partitions. In this paper we defined and investigated special spreads. Note that, for any odd prime power $q$, every special spread gives an equitable $(q^2+1)$-partition of $Sp(4,q)$, where each part is a coclique and there is a matching on the parts. As was discussed in this paper, it is possible to merge classes of this equitable partition to get many non-equivalent $s$-equitable 2-partitions of $Sp(4,q)$.  
It follows from \cite[Section 2.2.8, Case $d - e = s$]{BV22} (or it can be directly verified) that each of the parts $V_1$ and $V_2$ of an $s$-equitable partition in Theorem \ref{thm:DDG2} is an $h$-ovoid, where $h = \frac{q+1}{2}$. On the other hand, the complement of an $\frac{q+1}{2}$-ovoid in $Sp(4,q)$ is again an $\frac{q+1}{2}$-ovoid. Any pair consisting of an $\frac{q+1}{2}$-ovoid and its complement satisfies the conditions of Theorem \ref{thm:DDG2}. So the the following problem is of interest in context of this paper.

\begin{problem}
What are $\frac{q+1}{2}$-ovoids in $Sp(4,q)$? In particular, do the $\frac{q+1}{2}$-ovoids obtained from special spreads exhaust all $\frac{q+1}{2}$-ovoids in $Sp(4,q)$?    
\end{problem}

We did not search in the literature for counter examples to the last question in Problem 6, but the answer is certainly affirmative for $q=3$ as there is a unique 2-ovoid in $Sp(4,3)$, see e.g. Table 2 in \cite{BDB}.

A part of an $r$-equitable 2-partition of $Sp(2e,q)$ is known as a tight set (see \cite[Section 2.2.8, Section 2.5.9]{BV22}). In particular, each of the parts $V_1$ and $V_2$ of an $r$-equitable partition in Theorem \ref{thm:DDG4} is a $\frac{q^2+1}{2}$-tight set. (In general, the union of $i$ mutually disjoint cliques of size $\frac{q^e-1}{q-1}$ is always an $i$-tight set.) On the other hand, the complement of a $\frac{q^2+1}{2}$-tight set is again a $\frac{q^2+1}{2}$-tight set. Any pair consisting of a $\frac{q^2+1}{2}$-tight set and its complement satisfies the conditions of Theorem \ref{thm:DDG4}. A similar problem on tight sets in $Sp(4,q)$ is therefore also of interest in the context of this paper.

\begin{problem}
What are $\frac{q^2+1}{2}$-tight sets in $Sp(4,q)$? In particular, do the $\frac{q^2+1}{2}$-tight sets obtained from symplectic spreads exhaust all $\frac{q^2+1}{2}$-tight sets in $Sp(4,q)$?    
\end{problem}

\end{document}